\documentclass[twoside]{siamlatex}

\usepackage{bm}
\usepackage{graphicx}

\usepackage[bottom]{footmisc}
\usepackage{url}

\usepackage{times,amsbsy,amssymb,amsfonts,amscd,epic}
\usepackage{epsfig}
\usepackage{graphicx}
\usepackage{color}
\usepackage{geometry}
\usepackage{hyperref}

\usepackage[linesnumbered,ruled,vlined]{algorithm2e}
\usepackage{algorithmicx}
\usepackage{algpseudocode}
\usepackage{varwidth}
\usepackage{afterpage}

\newcommand{\V}[1]{\boldsymbol{#1}}

\newtheorem{thm}{Theorem}[section]
\newtheorem{lem}{Lemma}[section]

\renewcommand{\div}{{\rm div\,}}

\newcommand{\reff}[1]{(\ref{#1})}
\newcommand{\Proof}{\noindent {\bf Proof~}\ }
\newcommand{\Endproof}{$\hfill\Box$}

\newcommand{\be}{\begin{equation}}
\newcommand{\ee}{\end{equation}}
\newcommand{\ba}{\begin{array}}
\newcommand{\ea}{\end{array}}
\newcommand{\ben}{\begin{eqnarray}}
\newcommand{\een}{\end{eqnarray}}
\newcommand{\bn}{\begin{eqnarray*}}
\newcommand{\en}{\end{eqnarray*}}
\newcommand{\p}{\partial}
\newcommand{\bu}{{\bf u}}
\newcommand{\bv}{{\bf v}}
\newcommand{\bw}{{\bf w}}
\newcommand{\bbn}{{\bf n}}
\newcommand{\bX}{{\bf X}}
\newcommand{\ul}{\underline u}
\newcommand{\vl}{\underline v}

\newcommand{\Xl}{\underline X}
\newcommand{\cTfh}{{\cal T}_{f,h}}
\newcommand{\cTph}{{\cal T}_{p,h}}

\usepackage{graphicx}

\newtheorem{my assumption}{Assumption}

\title{Some Multilevel Decoupled Algorithms for a Mixed Navier-Stokes/Darcy Model}
\author{
Mingchao Cai \thanks{Department of Mathematics, Morgan State University, 1700 E Cold Spring Ln, Baltimore, MD 21251, USA.  \texttt{cmchao2005@gmail.com}. This work used the Extreme Science and Engineering Discovery Environment (XSEDE), which is supported by National Science Foundation grant number TG-DMS150025. This work has been submitted to Adv. Comput. Math. on Nov. 28, 2016.}
\and
Peiqi Huang \thanks{Department of Applied
Mathematics, Nanjing Forestry University, Nanjing 210037, People's  Republic of China.
\texttt{E-mail address: pqhuang1979@163.com}. This author's work is supported by the National Natural Science Foundation of China grants 11226309 and 11301267.}
\and
Mo Mu \thanks{Department of Mathematics, Hong Kong University of Science and Technology, Kowloon, Hong Kong. \texttt{mamu@ust.hk}. This work is supported in part by Hong Kong RGC Competitive Earmarked Research Grant HKUST603212.}
}
\begin{document}
\maketitle
\begin{center}
 \today
\end{center}
\begin{abstract}
In this work, several multilevel decoupled algorithms are proposed for a mixed Navier-Stokes/Darcy model. These algorithms are based on either successively or parallelly solving two linear subdomain problems after solving a coupled nonlinear coarse grid problem. Error estimates are given to demonstrate the approximation accuracy of the algorithms. Experiments based on both the first order and the second order discretizations are presented to show the effectiveness of the decoupled algorithms.
\end{abstract}

\begin{keywords}
Fluid flow coupled with porous media flow, Darcy law, Navier-Stokes equations, Interface coupling, Multilevel
algorithm, Decoupling, Linearization
\end{keywords}

\begin{AMS}
65F08, 65F10, 65N30, 65N55
\end{AMS}

\section{Introduction}
The coupling of incompressible fluid flow with porous media flow is an interesting but challenging topic. For describing the interactions of the fluid flow with the porous media flow, a coupled Stokes/Darcy or Naiver-Stokes/Darcy system is typically used as a macro-scale sharp interface model \cite{badia2009unified, badea2010numerical, burman2007unified, cai2009numerical, cai2009preconditioning, cao2010finite, chidyagwai2009solution, chidyagwai2011two, discacciati2002mathematical, discacciati2004convergence, discacciati2004domain, discacciati2007robin, ervin2009coupled, ervin2014approximation, girault2009dg, hou2016optimal, huang2012two, layton2003coupling, mu2007two, riviere2005locally, zhang2014two, zuo2014decoupling, zuo2015numerical}. The coupled Navier-Stokes/Darcy model is composed of a nonlinear Navier-Stokes equations for fluid flow, a Darcy law equation for porous media flow, plus certain interface conditions for describing the interactions of the different types of flows. Numerical methods for this model \cite{cai2009numerical, chidyagwai2009solution, girault2009dg, zuo2015numerical} usually result in a coupled and nonlinear saddle point problem, for which numerical difficulties increase as the mesh size decreases. 

\begin{figure}[htbp]
\centering
\includegraphics[scale=0.66]{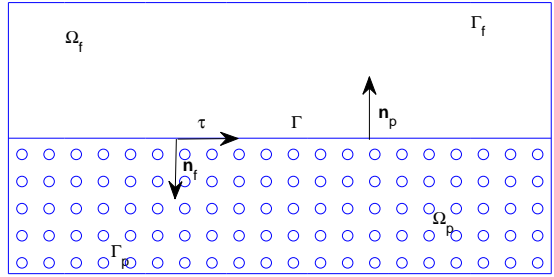}
\vspace{-2.0cm}
\caption{A global domain $\Omega$ consisting of a fluid region
$\Omega_f$ and a porous media region $\Omega_p$ separated by an
interface $\Gamma$.}
\label{domain}
\end{figure}

Let us consider a domain $\Omega \subset R^d$ ($d=2$ or 3), consisting of a fluid region $\Omega_f$ and a porous media region $\Omega_p$ separated by an interface $\Gamma$. As shown in Fig. \ref{domain}, $\Omega = \Omega_f \bigcup \Omega_p$ and $\Gamma=\overline{\Omega}_f \bigcap \overline{\Omega}_p $. The interface $\Gamma$ is assumed to be smooth enough \cite{girault2009dg}.

The fluid flow in $\Omega_f$ is governed by the steady state
Navier-Stokes equations:
\begin{equation}
\cases{ - \nu \Delta {\bf u}+\nabla p +\rho ({\bf u}\cdot
\nabla) {\bf u} = {\bf f} & $\forall {\bf x} \in \Omega_f$, \cr
\div {\bf u} = 0 & $\forall {\bf x} \in \Omega_f$, }
\label{NS_eqn}
\end{equation}
where $\rho$ is the density of the fluid flow, ${\bf u}$ is  the
velocity vector, $p$ is the pressure, ${\bf f}$ is the external
force, $\nu > 0$ is the viscosity coefficient.

In the porous media region $\Omega_p$, the governing equation becomes
\begin{equation}\label{Darcy}
\displaystyle
- \div \left( {\frac{{\bf K}}{n} \nabla \phi} \right) = f_p \quad \forall
{\bf x} \in \Omega_p.
\end{equation}
Here, $\phi$ is the piezometric head, $f_p$ is the source term due to injection or pump,
$n$ is the volumetric porosity, ${\bf K}$ is the hydraulic conductivity tensor of the porous media satisfying
$$
{\bf K}_{i j}= {\bf K}_{j i}, \quad \mbox{and} \quad \alpha_1 ({\bf
x}, {\bf x}) \le ({\bf K} {\bf x},{\bf x}) \le \alpha_2 ({\bf
x},{\bf x}) \quad \forall {\bf x} \in \Omega_p.
$$
Here, $\alpha_1$ and $\alpha_2$ are positive constants. Typically, ${\bf K}_{i j}$ is proportional to $\frac{\epsilon^2}{\nu}$ with $\epsilon$ being the characteristic length of the porous media. For simplicity, in this paper, we will assume that ${\bf K}=\frac{\epsilon^2}{\nu} {\bf I}$. In $\Omega_p$, the flow velocity and pressure can be calculated by
$$
{\bf u}_p= -\frac{\bf K}{n} \nabla \phi \quad \mbox{and}  \quad p_p= {\rho g} (\phi - z).
$$
Here, $z$, representing the elevation from a reference level, is assumed to be $0$, $p_p$ is the pressure in $\Omega_p$, and
$g$ is the gravity acceleration.

The key part of the coupled model is the transmission conditions at the interface, which describe the interaction mechanism of the two different types of flows. The following interface conditions have been extensively used and studied in the literature \cite{beavers1967boundary, discacciati2002mathematical, discacciati2004convergence, jager2000interface, layton2003coupling,
saffman1971boundary}:
              \be\label{Intrface-Condition}
               \left \{
                \ba {ll}
                 {\bf u}\cdot{\bf n}_f={\bf u}_p \cdot{\bf n}_f = -\frac{{\bf K}}{n} \nabla \phi \cdot{\bf n}_f,\\
                  -\nu(\nabla{\bf u}{\bf n}_f)\cdot{\bf n}_f+p=\rho g\phi,\\
                  -\nu(\nabla{\bf u}{\bf n}_f)\cdot{\V\tau}_i=\frac{\nu\alpha_{_{BJS}}}{\sqrt{\nu{\V\tau}_i\cdot {\bf K}{\V\tau}_i}}{\bf u}\cdot{\V\tau}_i, \quad i = 1, ..., d-1.
                \ea
               \right.
              \ee
Here, ${\bf n}_f$ is the unit outward normal directions on $\partial \Omega_f$ at $\Gamma$, $\{{\V\tau}_i\}_{i = 1}^{d-1}$ is the unit tangent vector on $\Gamma$, $\alpha_{_{BJS}}$ is a positive parameter depending on the properties of the porous medium. The first interface condition ensures the mass conservation across $\Gamma$. The second one is the balance of normal forces across the interface. The third condition is well known as Beavers-Joseph-Saffman's law \cite{beavers1967boundary, saffman1971boundary}, which states that the slip velocity is proportional to the shear stress along $\Gamma$.

For boundary conditions, without loss of generality, we impose homogeneous Dirichlet boundary conditions on $\Gamma_f = \partial \Omega_f / \Gamma$ and $\Gamma_p = \partial \Omega_p / \Gamma$:
\begin{equation}\label{Boundary-Condition}
\cases{ {\bf u} = {\bf 0} & on $\Gamma_f$, \cr
\phi = 0 & on $\Gamma_D$.
}
\end{equation}

The Finite Element method (FEM) discretization of the coupled Navier-Stokes/Darcy model will result in a coupled nonlinear saddle point problem, which is very difficult to solve. In this work, we are interested in developing decoupled and linearized methods so that they not only allow for easy and efficient implementation and software reuse, but also are numerically effective and efficient. We propose and investigate four multilevel decoupled algorithms. In all these algorithms, the coupled nonlinear system only needs to be solved on a very coarse grid level. After that, decoupled linearized Navier-Stokes and Darcy subproblems are solved on all the subsequently refined meshes. In {\it Algorithm A}, we solve a Darcy subproblem firstly and using the coarse grid solution to provide its boundary condition at the interface, and then solve a linearized Navier-Stokes problem using the Darcy problem to provide its boundary condition at the interface, and finally, on the same fine grid level, correct both the Darcy problem and the Navier-Stokes equations using the most updated subproblems to supplement the boundary conditions to each other at the interface. {\it Algorithm B} is similar to {\it Algorithm A}. Compared with {\it Algorithm A}, we only exchange the order of solving the linearized Navier-Stokes equations and the Darcy problem in {\it Algorithm B} \cite{zhang2014two}. In {\it Algorithm C}, on all fine grid levels, we use the previous level solution to provide boundary conditions for each subproblems and solve them in parallel \cite{mu2007two, cai2009numerical, zhang2014two, hou2016optimal}. In {\it Algorithm D}, on all fine grid levels, the correction step is only applied to the Navier-Stokes part, the boundary conditions of each subproblem are provided by using the most updated numerical solutions.

These multilevel algorithms are extended from the existing two-level algorithms \cite{huang2016newton, cai2009numerical, mu2007two, dai2008two, hou2016optimal, zuo2014decoupling, layton1993two, layton1995two, xu1996two}. However, the error estimates of the multilevel algorithms are much more difficult than those of two-level algorithms. In this paper, a theoretical analysis is given for {\it Algorithm A}. We apply mathematical induction method to give the estimates of the multilevel algorithm. Different from other existing papers, in which most of the researchers only analyze and test the first order discretization, our theory is valid not only for the first order discretization, but also valid for a general $k$-th order discretization. In particular, for both the first order and the second order discretizations, it is shown that if the mesh sizes of the two successive mesh levels are scaled with $h_l=h_{l-1}^2$, $l=1,2,\cdots,L$, then the energy norm errors in the final-step approximation are of optimal order. This means that the final approximation is of the same order of accuracy as the Finite Element approximation to $(\phi,\bu, p)$ obtained by solving exactly the coupled nonlinear system on the finest mesh. The results are similar to the so-called mesh independent principle justified for the multilevel algorithm for a single Navier-Stokes model by W. Layton \cite{layton1996multilevel, layton1998numerical}. The advantages of these multilevel algorithms are: they are numerically efficient because they enable the application of the most efficient and optimized local linear solvers on the fine grid that have been well developed for the linearized Navier-Stokes and Darcy models. Furthermore, in this work, we are interested in not only the mathematical analysis, but also the comparisons of different algorithms. Extensive numerical experiments for both the first order and the second order discretizations are provided to compare the different multilevel algorithms and to illustrate the effectiveness of these algorithms. In our numerical experiments, we firstly compare the algorithms in the two-level cases, then careful tests are designed to verify the theoretical predictions; some three-level experiments are also conducted to highlight the possible improvements of the theoretical analysis and the numerical algorithms.

The rest of the paper is organized as follows. The weak problem and a coupled and nonlinear algorithm are introduced in Section 2. Some multilevel algorithms are proposed in Section 3. Numerical analysis for {\it Algorithm A} is conducted in Section 4 to show that the decoupled and linearized multilevel algorithm retains the same order of approximation accuracy as the coupled and nonlinear algorithms if the scalings between the successive mesh levels are properly selected. In Section 5, we first compare the proposed two-level algorithms and then investigate the multilevel algorithms.

\section{Weak form and finite element approximations}

We begin with some notations. Let
\bn
{\bX}_f&=&\{{\bf v}\in {\bf H}^1(\Omega_f)=(H^1(\Omega_f))^d|\ {\bf v}={\bf 0}\ \mbox{on}\ \Gamma_f\}, \\
    Q&=&L^2(\Omega_f), \\
  X_p&=&\{\psi\in H^1(\Omega_p)|\ \psi=0\ \mbox{on}\ \Gamma_p\}
\en
be the functional spaces for $\bu$, $p$ and $\phi$, respectively. We denote $\Xl={\bX}_f \times X_p$. By multiplying test functions to \reff{NS_eqn} and \reff{Darcy}, integrating by parts and plugging in the interface boundary conditions \reff{Intrface-Condition}-\reff{Boundary-Condition}, the weak form of the coupled NS/Darcy model reads as: find $\ul=({\bf u},\phi)\in \Xl,p\in Q$ such that
              \be\label{Variational-Pro}
               \left \{
                \ba {ll}
                 a(\ul,\vl)+c(\bu,\bu,\bv)+b(\bv,p)=f(\vl)\qquad & \forall \vl=({\bf v},\psi)\in \Xl,\\
                 b(\bu,q)=0\qquad & \forall q\in Q,
                \ea
               \right.
              \ee
where
   \bn
     &a(\ul,\vl)=a_f({\bf u},{\bf v})+a_p(\phi,\psi)+a_\Gamma(\ul,\vl), \quad  b(\bv,p)=-\int_{\Omega_f}p\nabla\cdot{\bf v},\\
     &\emph{}c(\bu,\bv,{\bf w})=\rho\int_{\Omega_f}(\bu\cdot\nabla)\bv\cdot{\bf w}, \quad f(\vl)=\int_{\Omega_f}{\bf f}_f\cdot{\bf v}+\rho g \int_{\Omega_p}f_p\psi
   \en
with
\bn
a_f({\bf u},{\bf v})&=&\nu\int_{\Omega_f}\nabla\bu : \nabla\bv+ \sum_{i=1}^{d-1}\frac{ \nu\alpha_{_{BJS}}}{\sqrt{\nu{\V\tau}_i\cdot{\bf K}{\V\tau}_i}}\int_\Gamma({\bf u}\cdot{\V\tau}_i)({\bf v}\cdot{\V\tau}_i),\\
a_p(\phi,\psi)&=&\frac{\rho g}{n}\int_{\Omega_p}\nabla\psi\cdot{\bf K}\nabla\phi, ~~~~a_\Gamma(\ul,\vl)=\rho g\int_\Gamma(\phi{\bf v}-\psi{\bf u})\cdot{\bf n}_f.
\en
Here, $c(\bu,\bv,\bw)$ corresponds to the nonlinear term, $a_{f}({\bf u}, {\bf v}) + a_{p}(\phi,\psi)+a_\Gamma (u, v)$ is the corresponding bilinear form to the linear coupled Stokes/Darcy problem \cite{discacciati2002mathematical, mu2007two}. The following results have been well established: $a(\cdot,\cdot)$ is bounded and coercive; $b(\cdot, \cdot)$ is bounded and satisfies the inf-sup condition \cite{discacciati2002mathematical, discacciati2004domain, quarteroni1999domain}; and the nonlinear term satisfies the following estimates \cite{girault1986finite, layton1993two, layton1995two, layton1998two}.
\begin{lem}\label{Estimate_trilinear}
Suppose that the boundary of the domain $\Omega_f$ satisfies the strong Lipschitz condition of Adams \cite{adams1975sobolev}. We have
\bn
&&(a)\ |c(\bu,\bv,\bw)|\lesssim |\bu|_{1,\Omega_f}|\bv|_{1,\Omega_f}|\bw|_{1,\Omega_f}\qquad\qquad \forall \bu,\bv,\bw\in \bX_f,\\
&&(b)\ |c(\bu,\bv,\bw)|\lesssim |\bu|_{0,\Omega_f}|\bv|_{1,\Omega_f}\|\bw\|_{2,\Omega_f}\qquad\qquad \forall \bu,\bv\in \bX_f,\bw\in {\bf H}^2(\Omega_f),\\
&&(c)\ |c(\bu,\bv,\bw)|\lesssim |\bu|_{1,\Omega_f}|\bv|_{0,\Omega_f}\|\bw\|_{2,\Omega_f}\qquad\qquad \forall \bu,\bv\in \bX_f,\bw\in {\bf H}^2(\Omega_f).
\en
\end{lem}

Here and thereafter, we will use $a\lesssim b$ to denote that there exists a generic constant $C$, such that $a \le C b$. For the wellposedness of the coupled NS/Darcy model, we  refer to \cite{badea2010numerical, chidyagwai2009solution, discacciati2004domain, girault2009dg, zuo2015numerical}. It is shown that the coupled NS/Darcy problem (\ref{Variational-Pro}) is well-posed if the normal velocity across the interface is sufficiently small and the viscosity $\nu$ is sufficiently large. Moreover, there holds the a priori bound of the weak solution \cite{girault2009dg, zuo2015numerical}.

Now, we discuss the Finite Element approximations of problem (\ref{Variational-Pro}). For a subdomain $\Omega_d$, we denote $\|\cdot\|_{s,\Omega_d}$ and $|\cdot|_{s,\Omega_d}$ as the usual Sobolev norm and seminorm for $H^s(\Omega_d)$, respectively \cite{adams1975sobolev}. $(\cdot,\cdot)_{\Omega_d}$ represents the $L^2$ inner product on $\Omega_d$, where $\Omega_d$ can be the interface $\Gamma$ or one of the subdomains. We partition $\Omega_f$ and $\Omega_p$ by quasi-uniform regular triangulations $\cTfh$ and $\cTph$ with a characteristic meshsize $h$. Moreover, we assume that the two subdomain triangulations coincide at $\Gamma$. If a conventional conforming Finite Element method is applied to the model problem \reff{Variational-Pro}, the discrete problem reads as: Find $\ul_h=({\bf u}_h,\phi_h)\in \Xl_h={\bX}_{f,h}\times X_{p,h} \subset {\bX}_f\times X_p, p_h\in Q_h \subset Q$ such that
\be\label{Discrete-Pro}
               \left \{
                \ba {ll}
                 a(\ul_h,\vl_h)+c(\bu_h,\bu_h,\bv_h)+b(\bv_h,p_h)=f(\vl_h)\qquad & \forall \vl_h=({\bf v}_h,\psi_h)\in \Xl_h,\\
                 b(\bu_h,q_h)=0\qquad & \forall q_h\in Q_h.
                \ea
               \right.
\ee
Here, the FE pairs ${\bX}_{f,h}\times Q_h$ needs to be stable \cite{brezzi1991mixed, girault1986finite}, i.e., there exists a positive constant $\beta$ such that
\be\label{Discrete-LBB}
\sup_{\bv_h\in \bX_{f,h}}\frac{b(\bv_h,q_h)}{|\bv_h|_{1,\Omega_f}}\geq\beta\|q_h\|_{0,\Omega_f}\qquad \forall q_h\in Q_h.
\ee
We assume that the solution of (\ref{Variational-Pro}) is smooth enough and the FE spaces have the following typical approximation properties: let $k$ be a natural number, for all $(\bu,p)\in {\bf H}^{k+1}(\Omega_f)\cap\bX_f\times H^k(\Omega_f)$ and $\phi\in H^{k+1}(\Omega_p)\cap X_p$,
\be\label{FEMApproximation-up}
\inf_{\bv_h\in \bX_{f,h},q_h\in Q_h}\big\{h|\bu-\bv_h|_{1,\Omega_f}+\|\bu-\bv_h\|_{0,\Omega_f}+h\|p-q_h\|_{0,\Omega_f}\big\}
\lesssim h^{k+1}(|\bu|_{k+1,\Omega_f}+|p|_{k,\Omega_f});
\ee
\be\label{FEMApproximation-phi}
\inf_{\psi_h\in X_{p,h}}\big\{h|\phi-\psi_h|_{1,\Omega_p}+\|\phi-\psi_h\|_{0,\Omega_p}\big\}
\lesssim h^{k+1}|\phi|_{k+1,\Omega_p}.
\ee
There are several well-known Finite Element spaces satisfying the discrete {\it inf-sup} condition and the approximation properties (\ref{Discrete-LBB})-(\ref{FEMApproximation-phi}). For instance, if $k=1$, one can apply the Mini elements \cite{brezzi1991mixed, girault1986finite} in $\Omega_f$ and the piecewise linear elements in $\Omega_p$. If $k\ge 2$, the $k$-th order Taylor-Hood elements \cite{brezzi1991mixed, girault1986finite, taylor1973numerical} and $P_k$ elements can be applied in $\Omega_f$ and $\Omega_p$ respectively \cite{brezzi1991mixed, cai2008modeling, huang2016newton}. For simplicity, we will only consider the cases $k=1$ and $k=2$ for numerical experiments in this paper.

For the coupled discrete problem \reff{Discrete-Pro}, the energy norm error estimates can be derived by using a fixed-point framework \cite{cai2009numerical, girault1986finite}, the $L^2$ error analysis can be obtained by using the Aubin-Nitsche duality argument \cite{cai2009numerical}. In summary, we have
\begin{lem}\label{FEM-Err}
Let $(\bu,\phi,p)\in {\bf H}^{k+1}(\Omega_f)\times H^{k+1}(\Omega_p)\times H^k(\Omega_f)$ be the solution of the NS/Darcy model \reff{Variational-Pro} and $({\bf u}_h,\phi_h,p_h)$ be the Finite Element solution of \reff{Discrete-Pro}. Assuming that $\nu$ is sufficiently large and $h$ is sufficiently small, there holds the following energy norm estimate for the problem \reff{Discrete-Pro}.
\be\label{FEM-Err-H1}
  |\bu-\bu_h|_{1,\Omega_f}+|\phi-\phi_h|_{1,\Omega_p}+\|p-p_h\|_{0,\Omega_f}\lesssim h^k.
\ee
Moreover, we have the following $L^2$ error estimate:
\be\label{FEM-Err-L2}
  \|\bu-\bu_h\|_{0,\Omega_f}+\|\phi-\phi_h\|_{0,\Omega_p}\lesssim h^{k+1}.
\ee
\end{lem}
Furthermore, one can derive the a priori bound of the FE solution under the assumption that the viscosity $\nu$ is sufficiently large \cite{girault2009dg, zuo2015numerical}.

\section{Multilevel decoupled algorithms}
In this section, we introduce four multilevel decoupled algorithms for the coupled Navier-Stokes/Darcy model.
In the first step of all these algorithms, we solve the coupled nonlinear problem $(\ref{Discrete-Pro})$ on a coarse mesh level: find $u_H = ({\bf u}_H, \phi_H) \in {\Xl_H} \subset {\Xl}, ~ p_H \in Q_H \subset Q$ such that
\begin{equation}\label{couple_coarse_problem}
\cases{  a(\ul_H,\vl_H)+c(\bu_H,\bu_H,\bv_H) + b(v_H, p_H) = (f, v_H) & $\forall v_H =
({\bf v}_H, \psi_H) \in {\Xl_H}$, \cr b(u_H, q_H) = 0 & $\forall q_H \in
Q_H$.}
\end{equation}

In the following, for the ease of notations, we denote
$$
\tilde{a}_f (\bu, \bv, \bw) =a_f({\bf u},{\bf v})+c (\bu, \bv, \bw)+c (\bv, \bu, \bw).
$$

The first multi-level algorithm is actually an extension of the two-level algorithm developed in \cite{huang2016newton}. After solving the nonlinear coupled problem on a coarse grid level (cf. equation \reff{couple_coarse_problem}), the fine-level steps read as:

{\it Algorithm A}
\begin{algorithmic}[1]

\State Set $h_0=H$, $\phi^{h_0}_A=\phi_H$, $\bu^{h_0}_A=\bu_H$, and $p^{h_0}_A=p_H$.
\State For $l=1, 2,...,L$,

{\bf Step a}: Solve a Darcy problem on a fine grid: find $\phi^{*}_{A, h_l} \in X_{p,{h_l}} \supset X_{p, h_{l-1}}$ such that
\begin{equation}\label{FG_darcy1_A}
a_{p} (\phi^{*}_{A, h_l}, \psi_{h_l}) = (f_p, \psi_{h_l})  +\rho g(\psi_{h_l}, \bu^{h_{l-1}}_A\cdot\bbn_f)_{\Gamma}  \quad \forall \psi_{h_l} \in X_{p, h_l}.
\end{equation}

{\bf Step b}: Solve a linearized NS problem on a fine grid: find ${\bf u}^{*}_{A, h_l} \in {\bX}_{f, h_l}\supset {\bX}_{f, h_{l-1}}, ~
p^{*}_{A, h_l} \in Q_{h_l} \supset Q_{h_{l-1}}$ such that $\forall {\bf v}_{h_l} \in {\bX}_{f, h_l}$, $\forall q_{h_l} \in Q_{h_l}$,
\begin{equation}\label{FG_ns1_A}
\cases{\tilde{a}_f({\bf u}^{h_{l-1}}_{A}, {\bf u}^{*}_{A, h_l}, {\bf v}_{h_l}) + b({\bf v}_{h_l}, p^{*}_{A, h_l}) = (\tilde{\bf f}, {\bf v}_{h_l})-\rho g(\phi^*_{A, h_l},\bv_{h_l} \cdot\bbn_f)_{\Gamma}  \cr
b({\bf u}^{*}_{A, h_l}, q_{h_l}) = 0.
}
\end{equation}
Here, $(\tilde{\bf f}, {\bf v}_{h_l}) = ({\bf f}, {\bf v}_{h_l}) + c({\bf u}^{h_{l-1}}_{A},{\bf u}^{h_{l-1}}_{A},{\bf v}_{h_l})$.

{\bf Step c}: Correct the Darcy solution on the same fine grid: find $\phi^{h_l}_A \in X_{p,{h_l}}$ such that
\begin{equation}\label{FG_darcy2_A}
a_{p} (\phi^{h_l}_A, \psi_{h_l}) = (f_p, \psi_{h_l}) +\rho g(\psi_{h_l}, \bu^*_{A, h_l}\cdot\bbn_f)_{\Gamma}  \quad \forall \psi_{h_l} \in X_{p, h_l}.
\end{equation}

{\bf Step d}: Correct the NS solution on the same fine grid: find ${\bf u}^{h_l}_{A} \in {\bX}_{f, h_l}\supset {\bX}_{f, h_{l-1}}, p^{h_l}_{A} \in Q_{h_l} \supset Q_{h_{l-1}}$ such that $\forall {\bf v}_{h_l} \in {\bX}_{f, h_l}$, $\forall q_{h_l} \in Q_{h_l}$,
\begin{equation}\label{FG_ns2_A}
\cases{\tilde{a}_f({\bf u}^{h_{l-1}}_{A},{\bf u}^{h_l}_{A}, {\bf v}_{h_l}) + b({\bf v}_{h_l}, p^{h_l}_{A}) = (\bar{\bf f}, {\bf v}_{h_l}) -\rho g(\phi^{h_l}_A,\bv_{h_l} \cdot\bbn_f)_{\Gamma} \cr
b({\bf u}^{h_l}_{A}, q_{h_l}) = 0.
}
\end{equation}
~~~~Here,
$(\bar{\bf f}, {\bf v}_{h_l})=({\bf f}, {\bf v}_{h_l})+c(\bu^{h_{l-1}}_{A},\bu^*_{A, h_l},\bv_{h_l})+c(\bu^*_{A, h_l},\bu^{h_{l-1}}_A-\bu^*_{A, h_l},\bv_{h_l}).
$
\State End

\end{algorithmic}

In the second multi-level algorithm, different from {\it Algorithm A}, we exchange the order of solving the two subproblems on fine grid levels \cite{zuo2014decoupling}. Specifically, after solving the coupled nonlinear problem on a coarse grid level, the fine-level steps of the second multi-level algorithm read as:

{\it Algorithm B}
\begin{algorithmic}[1]

\State Set $h_0=H$, $\phi^{h_0}_B=\phi_H$, $\bu^{h_0}_B=\bu_H$, and $p^{h_0}_B=p_H$.
\State For $l=1, ..., L$.

{\bf Step a}: Solve a linearized NS problem on a fine grid: find ${\bf u}^{*}_{B,h_l} \in {\bX}_{f, h_l}\supset {\bX}_{f, h_{l-1}}, ~
p^{*}_{B, h_l} \in Q_{h_l} \supset Q_{h_{l-1}}$ such that
\begin{equation}\label{FG_ns1_B}
\cases{\tilde{a}_f({\bf u}^{ h_{l-1}}_{B},{\bf u}^{*}_{B, h_l}, {\bf v}_{h_l}) + b({\bf v}_{h_l}, p^{*}_{B, h_l}) = (\tilde{\bf f}, {\bf v}_{h_l})-\rho g(\phi_{B}^{h_{l-1}},\bv_{h_l} \cdot\bbn_f)_{\Gamma}  \cr
b({\bf u}^{*}_{B, h_l}, q_{h_l}) = 0. }
\end{equation}
Here, $(\tilde{\bf f}, {\bf v}_{h_l}) = ({\bf f}, {\bf v}_{h_l}) + c({\bf u}_{B}^{h_{l-1}},{\bf u}_{B}^{h_{l-1}},{\bf v}_{h_l})$.

{\bf Step b}: Solve a Darcy problem on a fine grid: find $\phi^{*}_{B, h_l} \in X_{p,{h_l}}$ such that
\begin{equation}\label{FG_darcy1_B}
a_{p} (\phi^{*}_{B, h_l}, \psi_{h_l}) = (f_p, \psi_{h_l})  +\rho g(\psi_{h_l}, \bu^*_{B, h_{j}}\cdot\bbn_f)_{\Gamma}  \quad \forall \psi_{h_l} \in X_{p, h_l}.
\end{equation}

{\bf Step c}: Correct the NS solution on the same fine grid: find ${\bf u}^{h_l}_B \in {\bX}_{f, h_l}\supset {\bX}_{f, h_{l-1}}, p^{h_l}_B \in Q_{h_l} \supset Q_{h_{l-1}}$ such that
\begin{equation}\label{FG_ns2_B}
\cases{\tilde{a}_f({\bf u}^{h_{l-1}}_B;{\bf u}^{h_l}_B, {\bf v}_{h_l}) + b({\bf v}_{h_l}, p^{h_l}_B) = (\bar{\bf f}, {\bf v}_{h_l}) -\rho g(\phi^*_{B, h_l}, \bv_{h_l} \cdot\bbn_f)_{\Gamma} \cr
b({\bf u}^{h_l}_B, q_{h_l}) = 0.
}
\end{equation}
Here, $(\bar{\bf f}, {\bf v}_{h_l})=({\bf f}, {\bf v}_{h_l})+c(\bu^{h_{l-1}}_B,\bu^*_{B, h_l},\bv_{h_l})+c(\bu^*_{B, h_l},\bu^{h_{l-1}}_B-\bu^*_{B, h_l},\bv_{h_l}). $

{\bf Step d}: Correct the Darcy solution on the same fine grid: find $\phi^{h_l}_B \in X_{p,{h_l}}$ such that
\begin{equation}\label{FG_darcy2_B}
a_{p} (\phi^{h_l}_B, \psi_{h_l}) = (f_p, \psi_{h_l}) +\rho g(\psi_{h_l}, \bu^{h_l}_{B}\cdot\bbn_f)_{\Gamma}  \quad \forall \psi_{h_l} \in X_{p, h_l}.
\end{equation}

\State End
\end{algorithmic}

In the third multilevel algorithm, after solving the coupled nonlinear problem on a coarse grid, we will solve the two subproblems in parallel on all fine grid levels. Specifically, the fine-level steps read as:

{\it Algorithm  C}
\begin{algorithmic}[1]

\State Set $h_0=H$, $\phi^{h_0}_C=\phi_H$, $\bu^{h_0}_C=\bu_H$, and $p^{h_0}_C=p_H$.
\State For $l=1, ..., L$,

Solve a linearized NS equation: find ${\bf u}^{h_l}_{C} \in {\bX}_{f, h_l}\supset {\bX}_{f, h_{l-1}}, ~p^{h_l}_{C} \in Q_{h_l} \supset Q_{h_{l-1}}$ such that
$\forall {\bf v}_{h_l} \in {\bX}_{f, h_l}$,  $\forall q_{h_l} \in Q_{h_l}$,
\begin{equation}\label{FG_nsC}
\cases{\tilde{a}_f({\bf u}^{h_{l-1}}_{C},{\bf u}^{h_l}_{C}, {\bf v}_{h_l})+ b({\bf v}_{h_l}, p^{h_l}_{C}) = ({\bf f}, {\bf v}_{h_l}) - \rho g (\phi^{h_{l-1}}_{C}, {\bf v}_{h_l}\cdot {\bbn}_f)_{\Gamma} \cr
b({\bf u}^{h_l}_{C}, q_{h_l}) = 0.
}
\end{equation}

Solve the local Darcy problem on a fine grid: find $\phi^{h_l}_C \in X_{p,{h_l}}$ such that
\begin{equation}\label{FG_darcyB2}
a_{p} (\phi^{h_l}_C, \psi_{h_l}) = (f_p, \psi_{h_l}) +  \rho g (\psi_{h_l}, {\bf u}^{h_{l-1}}_C \cdot {\bbn}_f)_{\Gamma} \quad \forall \psi_{h_l} \in X_{p, h_l}.
\end{equation}

\State End
\end{algorithmic}

In the last multi-level algorithm, we skip the correction step for the Darcy problem. After solving the coupled nonlinear problem on a coarse grid level, the fine-level steps of the algorithm reads as:

{\it Algorithm D}
\begin{algorithmic}[1]

\State Set $h_0=H$, $\phi^{h_0}_D=\phi_H$, $\bu^{h_0}_{D}=\bu_H$ and $p^{h_0}_D=p_H$.
\State For $l=1, ..., L$,

{\bf Step a}: Solve the Darcy problem on a fine grid: find $\phi^*_{D, h_l} \in X_{p,{h_l}} \supset X_{p, h_{l-1}}$ such that
\begin{equation}\label{FG_darcy1_D}
a_{p} (\phi^*_{D, h_l}, \psi_{h_l}) = (f_p, \psi_{h_l})  +\rho g(\psi_{h_l}, \bu^{h_{l-1}}_A\cdot\bbn_f)_{\Gamma}  \quad \forall \psi_{h_l} \in X_{p, h_l}.
\end{equation}

{\bf Step b}: Solve a linearized NS problem on a fine grid: find ${\bf u}^{*}_{D, h_l} \in {\bX}_{f, h_l}\supset {\bX}_{f, h_{l-1}}, ~
p^{*}_{D, h_l} \in Q_{h_l} \supset Q_{h_{l-1}}$ such that $\forall {\bf v}_{h_l} \in {\bX}_{f, h_l}$, $\forall q_{h_l} \in Q_{h_l}$,
\begin{equation}\label{FG_ns1_D}
\cases{\tilde{a}_f({\bf u}^{h_{l-1}}_{D}, {\bf u}^{*}_{D, h_l}, {\bf v}_{h_l}) + b({\bf v}_{h_l}, p^{*}_{D, h_l}) = (\tilde{\bf f}, {\bf v}_{h_l})-\rho g(\phi^*_{D, h_l}, \bv_{h_l} \cdot\bbn_f)_{\Gamma}  \cr
b({\bf u}^{*}_{D, h_l}, q_{h_l}) = 0. }
\end{equation}
Here, $(\tilde{\bf f}, {\bf v}_{h_l}) = ({\bf f}, {\bf v}_{h_l}) + c({\bf u}^{h_{l-1}}_{D},{\bf u}^{h_{l-1}}_{D},{\bf v}_{h_l})$.

{\bf Step c}: Correct the NS solution on the same fine grid: find ${\bf u}^{h_l}_{D} \in {\bX}_{f, h_l}\supset {\bX}_{f, h_{l-1}}, p^{h_l}_{D} \in Q_{h_l} \supset Q_{h_{l-1}}$ such that $\forall {\bf v}_{h_l} \in {\bX}_{f, h_l}$, $\forall q_{h_l} \in Q_{h_l}$,
\begin{equation}\label{FG_ns2_D}
\cases{\tilde{a}_f({\bf u}^{h_{l-1}}_{D},{\bf u}^{h_l}_{D}, {\bf v}_{h_l}) + b({\bf v}_{h_l}, p^{h_l}_{D}) = (\bar{\bf f}, {\bf v}_{h_l}) -\rho g(\phi^{h_{l-1}}_A,\bv_{h_l} \cdot\bbn_f)_{\Gamma} \cr
b({\bf u}^{h_l}_{D}, q_{h_l}) = 0.
}
\end{equation}
Here, $(\bar{\bf f}, {\bf v}_{h_l})=({\bf f}, {\bf v}_{h_l})+c(\bu^{h_{l-1}}_{D},\bu^*_{D, h_l},\bv_{h_l})+c(\bu^*_{D, h_l},\bu^{h_{l-1}}_A-\bu^*_{D, h_l},\bv_{h_l})$.

\State End
\end{algorithmic}

We see that when $L=1$, {\it Algorithm A} is reduced to the two-level algorithm developed in \cite{huang2016newton}, {\it Algorithm C} degenerates to the two-level algorithm proposed in \cite{mu2007two, cai2009numerical}. {\it Algorithm B} is an extension of the two-grid algorithm proposed in \cite{zhang2014two}. {\it Algorithm D} differs from {\it Algorithm A} in that there is no correction step for the Darcy problem. Intuitively, each of the above multilevel algorithms can be thought of as a recursive call of a certain two-level algorithm. Moreover, it is not difficult to see that {\it Algorithm A} and {\it Algorithm B} require more operation cost, while {\it Algorithm C} requires the least operation cost on every mesh level.

\section{Theoretical Analysis}

In this section, we only analyze the solution by the decoupled multilevel {\it Algorithm A}. For {\it Algorithm A} and {\it Algorithm B},  we will see that they produce almost the same accurate solution from our numerical experiments in Section 5. The analysis of {\it Algorithm C} in the linear case can be found in \cite{cai2012multilevel}.
As previously pointed out, when $L=1$, all the above algorithms degenerate to the two level algorithms. We firstly present the results for $L=1$ case, and then provide the error analysis  for analyzing the numerical solution on a general meshlevel $l$.

\subsection{Results for the two level algorithms}

For {\it Algorithm A} in the two-level case, we have the following results \cite{huang2016newton}.
\begin{lem}\label{T-lev-Err}
Let $H$ and $h$ be the coarse grid size and the fine grid size, i.e., $h_0=H$ and $h_1=h$, and let $(\phi,\bu,p)$, $(\phi^*_{A,h},\bu^*_{A,h},p^*_{A, h})$, and $(\phi^h_A,\bu^h_A,p^h_A)$ be defined by the problems \reff{Variational-Pro}, \reff{FG_darcy1_A}-\reff{FG_ns1_A}, and \reff{FG_darcy2_A}-\reff{FG_ns2_A}, respectively. Under the assumptions of Lemma \ref{FEM-Err}, there holds
     \be\label{T-lev-Err-phi-*}
       |\phi-\phi^*_{A, h}|_{1,\Omega_p} + |\bu-{\bf u}^*_{A,h}|_{1,\Omega_f}+ \|p-p^*_{A,h}\|_{0,\Omega_f}\lesssim H^{k+1}+h^k;
     \ee
\be\label{L2-estimates-*-0}
  \|{\bf u}-{\bf u}^*_{A, h}\|_{0,\Omega_f}\lesssim H^{2k+1}+H^{k+1}h+h^{k+1};
\ee
\be\label{T-lev-Err-phiup-h}
       |\phi-\phi^h_A|_{1,\Omega_p} +  |\bu-{\bf u}^h_A|_{1,\Omega_f}+  \|p-p^h_A\|_{0,\Omega_f} \lesssim H^{2k+1}+H^{k+1}h+h^k;
\ee
\be\label{L2-estimates-h-0}
  \|{\bf u}-{\bf u}^h\|_{0,\Omega_f}\lesssim H^{k+1}h^2+H^{k+1}h^k+H^{2(k+1)}+h^{k+1}.
\ee
\end{lem}

For {\it Algorithm C} in the two level case, the corresponding analysis for the linear case can be found in \cite{hou2016optimal}. In short, there holds
\be\label{T-lev-Err-C}
       |\phi-\phi^h_C|_{1,\Omega_p} +  |\bu-{\bf u}^h_C|_{1,\Omega_f}+  \|p-p^h_C\|_{0,\Omega_f} \lesssim H^{k+1} +h^k.
     \ee

{\bf Remarks}. For {\it Algorithm A}, we comment here that Lemma \ref{T-lev-Err} indicates that when $h=H^3$ if $k=1$ ($h=H^{5/2}$ if $k=2$), the final-step solution of {\it Algorithm A} possesses the same order accuracy as the Finite Element solution in the energy norm. In comparison, for {\it Algorithm C}, the theoretical estimates of energy norm errors in \reff{T-lev-Err-C} suggest that one needs to take the scaling $h=H^2$ if $k=1$ ($h=H^{3/2}$ if $k=2$). For {\it Algorithm A}, to ensure the final-step solutions have optimal $L^2$ norm errors, one has to take the scaling $h=H^2$ if $k=1, 2$; To ensure the intermediate-step solutions have optimal energy convergence, one has to take $h =H^\frac{k+1}{k}$ ($h=H^{3/2}$ if $k=2$); To ensure the intermediate-step solutions have optimal $L^2$ norm errors, the scaling between the two grid sizes has to be taken as $h =\max\{H^\frac{k+1}{k}, H^\frac{2k+1}{k+1} \}$ ($h=H^{3/2}$ if $k=1, 2$).

\subsection{Analysis of the multilevel algorithms}

The main purpose in this part is to show that the multilevel decoupled and linearized {\it Algorithm A}, with a properly chosen scalings of the two successive meshlevel sizes, is of the same order of approximation accuracy as the coupled and nonlinear algorithm. Note that {\it Algorithm A} may be viewed as an approximation to the coupled Finite Element algorithm, we will analyze
the difference between the solution by {\it Algorithm A} and the solution by using the nonlinear coupled algorithm. 


To estimate the $L^2$-error of the intermediate-step solution of {\it Algorithm A} on the $l$-th ($l\geq 1$) mesh level, we will consider the following the dual problem of the linearized problem: given ${\bf g}\in {\bf L}^2(\Omega_f)$, find $(\bw,r)\in\bX_f\times Q$ such that $\forall(\bv,q)\in\bX_f\times Q$
\ben\label{Duality-LinearizedNS}
  &&\tilde{a}_f (\bu, \bv, \bw)+b(\bv,r)+b(\bw,q)+\rho g(\phi-\phi^{*}_{A, h_l},\bw\cdot\bbn_f)_\Gamma=({\bf g},\bv)_{\Omega_f}.
\een
If the solution of the linearized coupled NS/Darcy model has the regularity $(\bu, \phi) \in(H^2(\Omega_f ))^d \times H^2(\Omega_p)$ as assumed in Lemma \ref{FEM-Err} and for $\nu$ sufficiently large, the two convection terms $c(\bu, \bv,\bw)$ and $c(\bv, \bu, \bw)$ in the linear dual problem (\ref{Duality-LinearizedNS}) can be properly bounded, and thus we may assume that the solution of the problem  (\ref{Duality-LinearizedNS}) is locally smooth and has the regularity
    \be\label{H2-Regularity}
      \|\bw\|_{2,\Omega_f}+\|r\|_{1,\Omega_f}\lesssim \|{\bf g}\|_{0,\Omega_f}.
    \ee

\begin{thm}\label{M-lev-Err-A}
Let $(\phi,\bu,p)$, $(\phi^*_{A, h_l},\bu^*_{A,h_l},p^*_{A, h_l})$ and $(\phi^{*}_{A, h_l},\bu^*_{A, h_l},p^*_{A, h_l},\phi^{h_l}_A, \bu^{h_l}_A,p^{h_l}_A)$ be defined by problem \reff{Variational-Pro} and \reff{FG_darcy1_A}-\reff{FG_ns2_A} (on a fine grid level with the grid size $h_l$), respectively. Under the assumptions of Lemma \ref{FEM-Err}, the following error estimates hold:
     \be\label{M-lev-Err-phi*}
       |\phi-\phi^{*}_{A, h_l}|_{1,\Omega_p}\lesssim h_{l-1} |{\bf u}_{h_l}-{\bu}^{h_{l-1}}_A|_{1,\Omega_f} + ||{\bf u}_{h_l}-{\bu}^{h_{l-1}}_A||_{0,\Omega_f} +h_l^k;
     \ee
     \ben\label{M-lev-Err-u*}
       |{\bf u}-{\bf u}^*_{A, h_l}|_{1,\Omega_f}+\|p-p^*_{A, h_l}\|_{0,\Omega_f} \lesssim |\phi_{h_l}-\phi^{*}_{A, h_l}|_{1,\Omega_p}+|\bu_{h_l} - \bu_A^{h_{l-1}}|^2_{1, \Omega_f} +h_l^k;
     \een
     \ben\label{L2-estimates-u*}
     \|\bu-\bu^*_{A, h_l}\|_{0,\Omega_f}  \lesssim  & h_l |\bu -\bu^*_{A, h_l}|_{1,\Omega_f} + |\bu - \bu^{h_{l-1}}_A|_{1, \Omega_f} ||\bu - \bu_A^{h_{l-1}}||_{0, \Omega_f}+h_l^{k+1},   
     \een
     \be\label{M-lev-Err-phi}
       |\phi-\phi^{h_l}_A|_{1,\Omega_p}\lesssim h_{l} |{\bf u}_{h_l}-{\bu}^*_{A, h_l}|_{1,\Omega_f} + ||{\bf u}_{h_l}-{\bu}^*_{A, h_l}||_{0,\Omega_f}+h_l^k;
     \ee
     \ben\label{M-lev-Err-u}
       |{\bf u}-{\bu}^{h_l}_A|_{1,\Omega_f}+\|p-p^{h_l}_A\|_{0,\Omega_f}\lesssim  |\phi_{h_l}-\phi^{h_l}_A|_{1,\Omega_p}
        \nonumber  \\
       + |\bu_{h_l} - \bu^*_{A, h_l}|_{1, \Omega_f}|\bu_{h_l}-\bu_A^{h_{l-1}}|_{1,\Omega_f}+h_l^k;
     \een
     \ben\label{L2-estimates-0}
      \|{\bf u}-{\bu}^{h_l}_A\|_{0,\Omega_f} & \lesssim & h_l \big( |\bu -\bu^{h_l}_{A}|_{1,\Omega_f}
      +|\bu-\bu^{h_{l-1}}_A|_{1,\Omega_f}|\bu-\bu^*_{A, h_l}|_{1,\Omega_f} \big) \nonumber \\
      && +\|\bu-\bu^{h_{l-1}}_A\|_{0,\Omega_f}|\bu-\bu^*_{A, h_l}|_{1,\Omega_f}+h_l^{k+1}.
     \een
\end{thm}


\Proof (i). The proof of \reff{M-lev-Err-phi*} is very similar to the estimate of $\phi$ in the two-grid algorithms developed in \cite{mu2007two, cai2009numerical}.
First, by taking $\vl_{h_l}=({\bf 0},\phi_{h_l}-\phi^*_{A, h_l})$ in \reff{Discrete-Pro} and comparing with the discrete model \reff{FG_darcy1_D}, we have
$$
a_p(\phi_{h_l}-\phi^*_{A, h_l}, \phi_{h_l} -\phi^*_{A, h_l})=\rho g\big(({\bf u}_{h_l}-{\bf u}^{h_{l-1}}_A)\cdot{\bf n}_f,\phi_{h_l}-\phi^*_{A, h_l}\big)_\Gamma.
$$
Let $\theta \in H^1(\Omega_f)$ be the solution of the problem:
              \bn
               \left \{
                \ba {ll}
                  -\Delta\theta=0 \quad & {\rm in} \quad\Omega_f,\\
                  \theta=\phi_{h_l}-\phi^*_{A, h_l} \quad & {\rm on} \quad\Gamma,\\
                  \theta=0 \quad & {\rm on} \quad\Gamma_f.
                \ea
               \right.
              \en
$\theta$ is the harmonic extension of $\phi_{h_l}-\phi^*_{A, h_l}$ to the fluid flow region and satisfies the following estimate \cite{mu2007two}.
$$
|\theta|_{1,\Omega_f}\lesssim\| \phi_{h_l}-\phi^*_{A, h_l} \|_{H^{1/2}_{00}(\Gamma)}\lesssim|\phi_{h_l}-\phi^*_{A, h_l}|_{1,\Omega_p}.
$$
Then, integrating by parts and noting that both $\bu_{h_l}$ and ${\bu}^*_{h_l}$ satisfy the discrete divergence-free property, we have,
for any $q_{h_{l-1}} \in Q_{h_{l-1}}$,
      \bn
      \big(({\bf u}_{h_l}-{\bu}^{h_{l-1}}_A)\cdot{\bf n}_f,\phi_{h_l}-\phi^*_{A, h_l}\big)_\Gamma &=&\big(\nabla\cdot({\bf u}_{h_l}-{\bu}^{h_{l-1}}_A),\theta- q_{h_{l-1}} \big)_{\Omega_f}+\big({\bf u}_{h_l}-{\bu}^{h_{l-1}}_A,\nabla\theta\big)_{\Omega_f}.
     \en
By applying the Cauchy-Schwarz inequality and the inequalities \reff{FEM-Err-L2}, there holds
     \bn
       |\phi_{h_l}-\phi^{*}_{A, h_l}|^2_{1,\Omega_p}
	   &\lesssim& a_p(\phi_{h_l}-\phi^{*}_{A, h_l},\phi_{h_l}-\phi^{*}_{A, h_l})\\
       &\lesssim&|{\bf u}_{h_l}-{\bu}^{h_{l-1}}_A|_{1,\Omega_f}\inf_{q_h\in Q_{h}}\|\theta-q_{h_{l-1}} \|_{0,\Omega_f}+\|{\bf u}_{h_l}-{\bu}^{h_{l-1}}_A\|_{0,\Omega_f}|\theta|_{1,\Omega_f}\\
       &\lesssim&\big( h_{l-1} | {\bf u}_{h_l}-{\bu}^{h_{l-1}}_A|_{1,\Omega_f}+\|{\bf u}_{h_l}-{\bu}^{h_{l-1}}_A\|_{0,\Omega_f} \big)|\theta|_{1,\Omega_f} \\
       &\lesssim&\big( h_{l-1} | {\bf u}_{h_l}-{\bu}^{h_{l-1}}_A|_{1,\Omega_f}+\|{\bf u}_{h_l}-{\bu}^{h_{l-1}}_A\|_{0,\Omega_f} \big)|\phi_{h_l}-\phi^{*}_{A, h_l}|_{1,\Omega_p}.
	   \en
By applying triangle inequality and the estimate of the Finite Element solution, we see that \reff{M-lev-Err-phi*} holds true.

(ii). We only provide a proof for the error estimate of ${\bf u}_{h_l}-{\bf u}^*_{A, h_l}$. Similar to the techniques used in \cite{mu2007two, cai2009numerical, zuo2014decoupling, zhang2014two}, the estimate for $p_{h_l} - p^*_{A, h_l}$ then follows from the discrete inf-sup condition and the estimate of ${\bf u}_{h_l}-{\bf u}^*_{A, h_l}$.
To prove (\ref{M-lev-Err-u*}), we compare the coupled nonlinear discrete problem \reff{Discrete-Pro} with the linearized Navier-Stokes model \reff{FG_ns1_A}. We see that
\ben\label{M-lev-Err-1}
   &&a_f({\bf u}_{h_l}-{\bf u}^*_{A, h_l},{\bf v}_{h_l})+c({\bf u}_{h_l}, {\bf u}_{h_l},{\bf v}_{h_l}) - [c({\bf u}^{h_{l-1}}_{A},  {\bf u}^*_{A, h_l}, {\bf v}_{h_l})
   +c({\bf u}^*_{A, h_l}, {\bf u}^{h_{l-1}}_A, {\bf v}_{h_l})   \nonumber \\
   &&-c({\bf u}_{A}^{h_{l-1}},  {\bf u}_A^{h_{l-1}}, {\bf v}_{h_l})]+b(\bv_{h_l},p_{h_l}-p^*_{A, h_l})=-\rho g(\phi_{h_l}-\phi^{h_l}_A,\bv_{h_l}\cdot\bbn_f)_\Gamma.
   \een
Taking ${\bf v}_{h_l}={\bf u}_{h_l}-{\bf u}^*_{A, h_l}$, due to the discrete divergence-free property of ${\bf u}_{h_l}$ and ${\bf u}^*_{A, h_l}$, there holds $b({\bf u}_{h_l}-{\bf u}^*_{A, h_l},p_{h_l}-p^*_{A, h_l})=0$. The interface term in \reff{M-lev-Err-1} can be controlled by $|\phi_{h_l}-\phi^{*}_{A, h_l}|_{1,\Omega_p} |{\bf u}_{h_l}-{\bf u}^*_{A, h_l}|_{1,\Omega_f}$. For the trilinear terms, it is easy to verify the following identity.
\begin{equation}\label{cf_identity}
\begin{array}{ccc}
~~c({\bf u}_{h_l}, {\bf u}_{h_l},{\bf u}_{h_l}-{\bf u}^*_{A, h_l}) - [c({\bf u}^{h_{l-1}}_{A},  {\bf u}^*_{A, h_l}, {\bf u}_{h_l}-{\bf u}^*_{A, h_l})+c({\bf u}^*_{A, h_l}, {\bf u}^{h_{l-1}}_A, {\bf u}_{h_l}-{\bf u}^*_{A, h_l}) \nonumber \\
-c({\bf u}_{A}^{h_{l-1}},  {\bf u}_A^{h_{l-1}}, {\bf u}_{h_l}-{\bf u}^*_{A, h_l})] = c({\bf u}_A^{h_{l-1}}, {\bf u}_{h_l}-{\bf u}^*_{A, h_l}, {\bf u}_{h_l}-{\bf u}^*_{A, h_l}) \nonumber \\
+c({\bf u}_{h_l}-{\bf u}^*_{A, h_l}, {\bf u}_{A}^{h_{l-1}}, {\bf u}_{h_l}-{\bf u}^*_{A, h_l})+c({\bf u}_{h_l}-{\bf u}_A^{h_{l-1}}, {\bf u}_{h_l}-{\bf u}_A^{h_{l-1}}, {\bf u}_{h_l}-{\bf u}^*_{A, h_l}).
\end{array}
\end{equation}
Note that the viscosity $\nu$ is sufficiently large, then the weak solution, the Finite Element solution as well as the multilevel solution have the a-priori bounds \cite{girault2009dg, zuo2015numerical}. Then, roughly speaking, the following inequality holds true.
$$
|c({\bu}^{h_{l-1}}_A,{\bf u}-{\bf u}^*_{A, h_l},{\bf u}-{\bf u}^*_{A, h_l})+c({\bf u}-{\bf u}^*_{A, h_l},{\bu}^{h_{l-1}}_A,{\bf u}-{\bf u}^*_{A, h_l}) | \le \frac{\nu}{2}|{\bf u}_{h_l}-{\bf u}^*_{A, h_l}|^2_{1,\Omega_f}.
$$
Thus, by using  \reff{M-lev-Err-1}, \reff{cf_identity}, and Lemma 2.1, we see that
\bn
\frac{\nu}{2} |{\bf u}_{h_l}-{\bf u}^*_{A, h_l}|^2_{1,\Omega_f} &\le &a_f({\bf u}_{h_l}-{\bf u}^*_{A, h_l},{\bf u}_{h_l}-{\bf u}^*_{A, h_l}) +c({\bf u}_A^{h_{l-1}}, {\bf u}_{h_l}-{\bf u}^*_{A, h_l}, {\bf u}_{h_l}-{\bf u}^*_{A, h_l}) \nonumber \\
&& +c({\bf u}_{h_l}-{\bf u}^*_{A, h_l}, {\bf u}_{A}^{h_{l-1}}, {\bf u}_{h_l}-{\bf u}^*_{A, h_l})  \\
& = &  - \rho g (\phi_{h_l}-\phi^*_{A, h_l}, ({\bf u}_{h_l}-{\bf u}^*_{A, h_l}) \cdot {\bf n}_f)_{\Gamma} -c({\bf u}_{h_l}-{\bf u}_A^{h_{l-1}}, {\bf u}_{h_l}-{\bf u}_A^{h_{l-1}}, {\bf u}_{h_l}-{\bf u}^*_{A, h_l}) \\
&\lesssim & |\phi_{h_l}-\phi^{*}_{A, h_l}|_{1,\Omega_p} |{\bf u}_{h_l}-{\bf u}^*_{A, h_l}|_{1,\Omega_f}+|\bu_{h_l} - \bu_A^{h_{l-1}}|^2_{1, \Omega_f}  |{\bf u}_{h_l}-{\bf u}^*_{A, h_l}|_{1,\Omega_f}.
\en
Hence, by applying the triangle inequality, $|{\bf u}-{\bf u}^*_{A, h_l}|_{1,\Omega_f} \le |{\bf u} - {\bf u}_{h_l}|_{1,\Omega_f} + |{\bf u}_{h_l}-{\bf u}^*_{A, h_l}|_{1,\Omega_f}$, and the energy norm estimate of the Finite Element solution \reff{FEM-Err-H1}, we see that the inequality \reff{M-lev-Err-u*} hods true.

(iii). For estimating the $L^2$ error of the intermediate-step solution, we set ${\bf g}=\bu-\bu^*_{A, h_l}$ and $(\bv,q)=(\bu-\bu^*_{A, h_l},p-p^*_{A, h_l})$ in \reff{Duality-LinearizedNS}, and then splitting the two trilinear terms into four terms, we obtain
\ben\label{L2-estimates-1}
  \|\bu-\bu^*_{A, h_l}\|^2_{0,\Omega_f}&=&a_f(\bu-\bu^*_{A, h_l},\bw)+c(\bu^{h_{l-1}}_A,\bu-\bu^*_{A, h_l},\bw)+c(\bu-\bu^*_{A, h_l},\bu^{h_{l-1}}_A,\bw)\nonumber\\
   &&+b(\bu-\bu^*_{A, h_l},r)+b(\bw,p-p^*_{A, h_l})+\rho g(\phi-\phi^*_{A, h_l},\bw\cdot\bbn_f)_\Gamma\nonumber\\
   &&+c(\bu-\bu^{h_{l-1}}_A,\bu-\bu^*_{A, h_l},\bw)+c(\bu-\bu^*_{A, h_l},\bu-\bu^{h_{l-1}}_A,\bw).
\een
Taking $\vl=({\bf v}_{h_l},0)$ in \reff{Variational-Pro} and subtracting with \reff{FG_ns1_A}, we obtain
   \ben\label{L2-estimates-2}
   &&a_f({\bf u}-{\bf u}^*_{A, h_l},{\bf v}_{h_l})+c({\bu}^{h_{l-1}}_A,{\bf u}-{\bf u}^*_{A, h_l},{\bf v}_{h_l})+c({\bf u}-{\bf u}^*_{A, h_l},{\bu}^{h_{l-1}}_A,{\bf v}_{h_l})+b(\bv_{h_l},p-p^*_{A, h_l})\nonumber\\
   &&\hspace{0.3cm}~~~~~~~~+b({\bf u}-{\bf u}^*_{A, h_l},q_{h_l})=-\rho g(\phi-\phi^*_{A, h_l},\bv_{h_l}\cdot\bbn_f)_\Gamma-c({\bf u}-{\bu}^{h_{l-1}}_A,{\bf u}-{\bu}^{h_{l-1}}_A,{\bf v}_{h_l}).
   \een
Subtracting \reff{L2-estimates-2} from \reff{L2-estimates-1}, we have
\bn
  \|\bu-\bu^*_{A, h_l}\|^2_{0,\Omega_f}&=&a_f(\bu-\bu^*_{A, h_l},\bw-{\bf v}_{h_l})+c(\bu^{h_{l-1}}_A,\bu-\bu^*_{A, h_l},\bw-{\bf v}_{h_l})+c(\bu-\bu^*_{A, h_l},\bu^{h_{l-1}}_A,\bw-{\bf v}_{h_l})\\
   &&+b(\bu-\bu^*_{A, h_l},r-q_{h_l})+b(\bw-{\bf v}_{h_l},p-p^*_{A, h_l})+\rho g(\phi-\phi^*_{A, h_l},(\bw-{\bf v}_{h_l})\cdot\bbn_f)_\Gamma\\
   &&+c(\bu-\bu^{h_{l-1}}_A,\bu-\bu^*_{A, h_l},\bw)+c(\bu-\bu^*_{A, h_l},\bu-\bu^{h_{l-1}}_A,\bw)\\
   &&+c({\bf u}-{\bu}^{h_{l-1}}_A,{\bf u}-{\bu}^{h_{l-1}}_A,\bw-{\bf v}_{h_l})-c({\bf u}-{\bu}^{h_{l-1}}_A,{\bf u}-{\bu}^{h_{l-1}}_A,\bw)\\
   &\lesssim&|\bu-\bu^*_{A, h_l}|_{1,\Omega_f}|\bw-\bv_{h_l}|_{1,\Omega_f}+2|\bu^{h_{l-1}}_A|_{1,\Omega_f}|\bu-\bu^*_{A, h_l}|_{1,\Omega_f}|\bw-{\bf v}_{h_l}|_{1,\Omega_f}\\
   &&+|\bu-\bu^*_{A, h_l}|_{1,\Omega_f}\|r-q_{h_l}\|_{0,\Omega_f}+|\bw-\bv_{h_l}|_{1,\Omega_f}\|p-p^*_{A, h_l}\|_{0,\Omega_f}\\
   &&+|\phi-\phi^*_{A, h_l}|_{1,\Omega_p}|\bw-\bv_{h_l}|_{1,\Omega_f}
   +2|\bu-\bu^*_{A, h_l}|_{1,\Omega_f}\|\bu-\bu^{h_{l-1}}_A\|_{0,\Omega_f}\|\bw\|_{2,\Omega_f}\\
   &&+|\bu-\bu^{h_{l-1}}_A|^2_{1,\Omega_f}|\bw-\bv_{h_l}|_{1,\Omega_f}
   +|\bu-\bu^{h_{l-1}}_A|_{1,\Omega_f}\|\bu-\bu^{h_{l-1}}_A\|_{0,\Omega_f}\|\bw\|_{2,\Omega_f} \\
   & \lesssim & (|\phi -\phi^*_{A, h_l} |_{1,\Omega_f}+|\bu -\bu^*_{A, h_l}|_{1,\Omega_f}+\|p - p^*_{A, h_l}\|_{0,\Omega_f}) |\bw-{\bf v}_{h_l}|_{1, \Omega_f} \\
   && +(|\bu - \bu^{h_{l-1}}_A|_{1, \Omega_f}+|\bu - \bu^*_{A, h_l}|_{1, \Omega_f} )||\bu - \bu_A^{h_{l-1}}||_{0, \Omega_f}  ||\bw||_{2, \Omega_f} \\
   && + |\bu-\bu^*_{A, h_l}|_{1,\Omega_f}\|r-q_{h_l}\|_{0,\Omega_f}.
\en
Here, in the last inequality we have used the Cauchy-Schwarz inequality and the estimate in Lemma \ref{Estimate_trilinear} for the trilinear term, we have dropped those higher order terms as $k \ge 1$ and $h_l \le h_{l-1}$. By the approximation error estimate \reff{FEMApproximation-up}-\reff{FEMApproximation-phi}, discarding the terms which are of the same order or higher order errors (for example, $|\phi -\phi^*_{A, h_l} |_{1,\Omega_f}$ and $\|p - p^*_{A, h_l}\|_{0,\Omega_f}$ are of the same order as $|\bu -\bu^*_{A, h_l}|_{1,\Omega_f}$, and $|\bu - \bu^*_{A, h_l}|_{1, \Omega_f} \le |\bu - \bu^{h_{l-1}}_A|_{1, \Omega_f}$), using Lemma \ref{FEM-Err} and the estimate \reff{H2-Regularity}, it follows that the $L^2$ error estimate for $\bu^*_{A, h_l}$ holds true.

(iv). The estimate of (\ref{M-lev-Err-phi}) is similar to the estimate of (\ref{M-lev-Err-phi*}). Taking $\vl_h=({\bf 0},\phi_{h_l}-\phi^{h_l}_A)$ in (\ref{Discrete-Pro}) and
comparing with the discrete model \reff{FG_darcy2_A}, we have
$$
a_p(\phi_{h_l}-\phi^{h_l}_A,\phi_{h_l}-\phi^{h_l}_A)=\rho g\big(({\bf u}_{h_l}-{\bu}^*_{A, h_l})\cdot{\bf n}_f,\phi_{h_l}-\phi^{h_l}_A\big)_\Gamma.
$$
Let $\theta \in H^1(\Omega_f)$ be a harmonic extension of $\phi_{h_l}-\phi^{h_l}_A$ to the fluid flow region with the Dirichlet data at $\Gamma$ being equal to equal to $\phi_{h_l}-\phi^{h_l}_A$.
Then, we have
   $$|\theta|_{1,\Omega_f}\lesssim\|\phi_{h_l}-\phi^{h_l}_A\|_{H^{1/2}_{00}(\Gamma)}\lesssim|\phi_{h_l}-\phi^{h_l}_A|_{1,\Omega_p},$$
Note that for any $q_{h_l}\in Q_{h_l}$, there holds
     \bn
       \big(({\bf u}_{h_l}-{\bu}^*_{A, h_l})\cdot{\bf n}_f,\phi_{h_l}-\phi^{h_l}_A\big)_\Gamma&=&\big(({\bf u}_{h_l}-{\bu}^*_{A, h_l})\cdot{\bf n}_f,\theta\big)_{\p\Omega_f}\\
       &=&\big(\nabla\cdot({\bf u}_{h_l}-{\bu}^*_{A, h_l}),\theta-q_{h_l}\big)_{\Omega_f}+\big({\bf u}_{h_l}-{\bu}^*_{A, h_l},\nabla\theta\big)_{\Omega_f}.
     \en
Here, in the last equality, we have used the discrete divergence-free property for ${\bf u}_{h_l}$ and ${\bu}^*_{A, h_l}$. Therefore,
     \begin{equation}\label{est_phi_A}
     \begin{array}{lll}
       |\phi_{h_l}-\phi^{h_l}_A|^2_{1,\Omega_p}&\lesssim&a_p(\phi_{h_l}-\phi^{h_l}_A,\phi_{h_l}-\phi^{h_l}_A)\\
       &\lesssim&| {\bf u}_{h_l}-{\bu}^*_{A, h_l} |_{1,\Omega_f}\inf_{q_{h_l}\in Q_{h_l}}\|\theta-q_{h_l}\|_{0,\Omega_f}+\|{\bf u}_{h_l}-{\bu}^*_{A, h_l}\|_{0,\Omega_f}|\theta|_{1,\Omega_f}\\
       &\lesssim&\big(h_l|{\bf u}_{h_l}-{\bu}^*_{A, h_l}|_{1,\Omega_f}+\|{\bf u}_{h_l}-{\bu}^*_{A, h_l}\|_{0,\Omega_f}\big)|\theta|_{1,\Omega_f}\\
     \end{array}
     \end{equation}
We see that the estimate \reff{M-lev-Err-phi} holds true.

(v). Now, we estimate the error of ${\bf u}_{h_l}-{\bu}^{h_l}_A$ in the energy norm. Taking $\vl_{h_l}=({\bf v}_{h_l},0)$ in \reff{Discrete-Pro} and comparing with \reff{FG_ns2_A} on the $l$-th level mesh, and splitting the trilinear terms in the right hand side, we obtain
   \ben\label{M-lev-Err-4}
   &&a_f({\bf u}_{h_l}-{\bu}^{h_l}_A,{\bf v}_{h_l})+c({\bu}^{h_{l-1}}_A,{\bf u}_{h_l}-{\bu}^{h_l}_A,{\bf v}_{h_l})+c({\bf u}_{h_l}-{\bu}^{h_l}_A,{\bu}^{h_{l-1}}_A,{\bf v}_{h_l}) +b({\bf v}_{h_l},p_{h_l}-p^{h_l}_A) \nonumber\\
   &&\hspace{1cm}=-\rho g(\phi_{h_l}-\phi^{h_l}_A,\bv_{h_l}\cdot\bbn_f)_\Gamma+c({\bf u}_{h_l}-{\bu}^{h_{l-1}}_A,{\bf u}^*_{A, h_l}-{\bf u}_{h_l},{\bf v}_{h_l})\nonumber\\
   &&\hspace{2cm}+c({\bf u}_{h_l}-{\bf u}^*_{A, h_l},{\bu}^{h_{l-1}}_A-{\bf u}_{h_l},{\bf v}_{h_l})+c({\bf u}^*_{A, h_l}-{\bf u}_{h_l},{\bf u}^*_{A, h_l}-{\bf u}_{h_l},{\bf v}_{h_l}).
   \een
Similar to the proof of \reff{M-lev-Err-u*}, letting ${\bf v}_{h_l}={\bf u}_{h_l}-{\bu}^{h_l}_A$, we have
   \ben\label{M-lev-Err-5}
    \frac{\nu}{2} |{\bf u}_{h_l}-{\bu}^{h_l}_A|^2_{1,\Omega_f} &\le& a_f({\bf u}_{h_l}-{\bu}^{h_l}_A,{\bf u}_{h_l}-{\bu}^{h_l}_A) +c({\bf u}_A^{h_{l-1}}, {\bf u}_{h_l}-{\bu}^{h_l}_A, {\bf u}_{h_l}-{\bu}^{h_l}_A)  \nonumber \\
    &&  +c({\bf u}_{h_l}-{\bu}^{h_l}_A,{\bu}^{h_{l-1}}_A, {\bf u}_{h_l}-{\bu}^{h_l}_A) \nonumber \\
     &=&-\rho g(\phi_{h_l}-\phi^{h_l}_A,({\bf u}_{h_l}-{\bu}^{h_l}_A)\cdot\bbn_f)_\Gamma- c({\bf u}_{h_l}-{\bu}^{h_{l-1}}_A,{\bf u}^*_{A, h_l}-{\bf u}_{h_l},{\bf u}_{h_l}-{\bu}^{h_l}_A)\nonumber\\
      &&+c({\bf u}_{h_l}-{\bf u}^*_{A, h_l},{\bu}^{h_{l-1}}_A-{\bf u}_{h_l},{\bf u}_{h_l}-{\bu}^{h_l}_A)+c({\bf u}^*_{A, h_l}-{\bf u}_{h_l},{\bf u}^*_{A, h_l}-{\bf u}_{h_l},{\bf u}_{h_l}-{\bu}^{h_l}_A)
      \een
The right hand side of \reff{M-lev-Err-5} is bounded by
\ben \label{M-lev-Err-final}
&&||\phi_{h_l}-\phi^{h_l}_A\|_{0,\Gamma} ||{\bf u}_{h_l}-{\bu}^{h_l}_A\|_{0,\Gamma}+|{\bf u}_{h_l}-{\bf u}^*_{A, h_l}|^2_{1,\Omega_f}|{\bf u}_{h_l}-{\bu}^{h_l}_A|_{1,\Omega_f}  \nonumber\\
&&+2|{\bf u}_{h_l}-{\bf u}^*_{A, h_l}|_{1,\Omega_f}|{\bf u}_{h_l}-{\bu}^{h_{l-1}}_A|_{1,\Omega_f}|{\bf u}_{h_l}-{\bu}^{h_l}_A|_{1,\Omega_f}  \nonumber \\
&\lesssim&\big(|\phi_{h_l}-\phi^{h_l}_A|_{1,\Omega_p} +|{\bf u}_{h_l}-{\bf u}^*_{A, h_l}|^2_{1,\Omega_f}+|{\bf u}_{h_l}-{\bf u}^*_{A, h_l}|_{1,\Omega_f}|{\bf u}_{h_l}-{\bu}^{h_{l-1}}_A|_{1,\Omega_f} \big)|{\bf u}_{h_l}-{\bu}^{h_l}_A|_{1,\Omega_f}. \nonumber
\een
Applying the triangle inequality, the energy norm error estimate of Finite Element solution (cf. Lemma \ref{FEM-Err}), then discarding the terms which are of the same order or higher order errors in \reff{M-lev-Err-4} and \reff{M-lev-Err-5}, we see that \reff{M-lev-Err-u} holds true.

(vi). For estimating the $L^2$-error of $\bu-\bu^{h_l}_A$, let us consider the dual problem of a linearized problem, which is similar to problem (\ref{Duality-LinearizedNS}).
\ben\label{Duality-LinearNS}
  &&\tilde{a}_f (\bu, \bv, \bw)+b(\bv,r)+b(\bw,q)+\rho g(\phi-\phi^{h_l}_A,\bw\cdot\bbn_f)_\Gamma=({\bf g},\bv)_{\Omega_f}.
\een
Moreover, we assume that a regularity estimate which is similar to \reff{H2-Regularity} holds. Setting ${\bf g}=\bu-\bu^{h_l}_A$ and $(\bv,q)=(\bu-\bu^{h_l}_A,p-p^{h_l}_A)$ in \reff{Duality-LinearNS}, splitting the two trilinear terms in \reff{Duality-LinearNS} into four terms,
we have
\ben\label{L2-estimates-4}
  \|\bu-\bu^{h_l}_A\|^2_{0,\Omega_f}&=&a_f(\bu-\bu^{h_l}_A,\bw)+c(\bu^{h_{l-1}}_A,\bu-\bu^{h_l}_A,\bw)+c(\bu-\bu^{h_l}_A,\bu^{h_{l-1}}_A,\bw)\nonumber\\
   &&+b(\bu-\bu^{h_l}_A,r)+b(\bw,p-p^{h_l}_A)+\rho g(\phi-\phi^{h_l}_A,\bw\cdot\bbn_f)_\Gamma\nonumber\\
   &&+c(\bu-\bu^{h_{l-1}}_A,\bu-\bu^{h_l}_A,\bw)+c(\bu-\bu^{h_l}_A,\bu-\bu^{h_{l-1}}_A,\bw).
\een
Taking $\vl=({\bf v}_{h_l},0)$ in \reff{Variational-Pro} and subtracting with \reff{FG_ns2_A}, we obtain
   \ben\label{L2-estimates-5}
   &&a_f({\bf u}-{\bu}^{h_l}_A,{\bf v}_{h_l})+c({\bu}^{h_{l-1}}_A,{\bf u}-{\bu}^{h_l}_A,{\bf v}_{h_l})+c({\bf u}-{\bu}^{h_l}_A,{\bu}^{h_{l-1}}_A,{\bf v}_{h_l})+b(\bv_{h_l},p-p^{h_l}_A)\nonumber\\
   &&\hspace{0.3cm}+b({\bf u}-{\bu}^{h_l}_A,q_{h_l})=-\rho g(\phi-\phi^{h_l}_A,\bv_{h_l}\cdot\bbn_f)_\Gamma-c({\bf u}-{\bu}^{h_{l-1}}_A,{\bf u}-{\bf u}^*_{A, h_l},{\bf v}_{h_l})\nonumber\\
   &&\hspace{2cm}-c({\bf u}-{\bf u}^*_{A, h_l},{\bf u}-{\bu}^{h_{l-1}}_A,{\bf v}_{h_l})+c(\bu-{\bf u}^*_{A, h_l},\bu-{\bf u}^*_{A, h_l},{\bf v}_{h_l}).
   \een
Combining \reff{L2-estimates-4} and \reff{L2-estimates-5}, we arrive at
\bn
  \|\bu-\bu^{h_l}_A\|^2_{0,\Omega_f}&=&a_f(\bu-\bu^{h_l}_A,\bw-{\bf v}_{h_l})+c(\bu^{h_{l-1}}_A,\bu-\bu^{h_l}_A,\bw-{\bf v}_{h_l})+c(\bu-\bu^{h_l}_A,\bu^{h_{l-1}}_A,\bw-{\bf v}_{h_l})\\
   &&+b(\bu-\bu^{h_l}_A,r-q_{h_l})+b(\bw-{\bf v}_{h_l},p-p^{h_l}_A)+\rho g(\phi-\phi^{h_l}_A,(\bw-{\bf v}_{h_l})\cdot\bbn_f)_\Gamma\\
   &&+c(\bu-\bu^{h_{l-1}}_A,\bu-\bu^{h_l}_A,\bw)+c(\bu-\bu^{h_l}_A,\bu-\bu^{h_{l-1}}_A,\bw)\\
   &&-c(\bu-\bu^{h_{l-1}}_A,\bu-\bu^*_{A, h_l},\bv_{h_l})-c(\bu-\bu^*_{A, h_l},\bu-\bu^{h_{l-1}}_A,\bv_{h_l})\\
   &&+c({\bf u}-{\bf u}^*_{A, h_l},{\bf u}-{\bf u}^*_{A, h_l},{\bf v}_{h_l}).
\en
Similar to the proof for $\|\bu-\bu^*_{A, h_l}\|^2_{0,\Omega_f}$, there holds
\bn\label{L2-estimates-6}
  \|\bu-\bu^{h_l}_A\|^2_{0,\Omega_f}  &\lesssim & \big( |\phi-\phi^{h_l}_A|_{1,\Omega_p}+|\bu-\bu^{h_l}_A|_{1,\Omega_f} + \|p-p^{h_l}_A\|_{0,\Omega_f} \big) |\bw-\bv_{h_l}|_{1,\Omega_f} +|\bu-\bu^{h_l}_A|_{1,\Omega_f}\|r-q_{h_l}\|_{0,\Omega_f} \\
   && +2\|\bu-\bu^{h_{l-1}}_A\|_{0,\Omega_f} |\bu-\bu^{h_l}_A |_{1,\Omega_f} \|\bw\|_{2,\Omega_f}\\
   &&+c({\bf u}-{\bf u}^*_{A, h_l},{\bf u}-{\bf u}^*_{A, h_l},{\bf v}_{h_l})-c(\bu-\bu^{h_{l-1}}_A,\bu-\bu^*_{A, h_l},\bv_{h_l})-c(\bu-\bu^*_{A, h_l},\bu-\bu^{h_{l-1}}_A,\bv_{h_l}). \nonumber
\en
Then, the estimates for the last three terms in \reff{L2-estimates-6} are:
\bn
|c({\bf u}-{\bf u}^*_{A, h_l},{\bf u}-{\bf u}^*_{A, h_l},{\bf v}_{h_l})|
  &\leq&|c({\bf u}-{\bf u}^*_{A, h_l},{\bf u}-{\bf u}^*_{A, h_l},\bw-{\bf v}_{h_l})|
    +|c({\bf u}-{\bf u}^*_{A, h_l},{\bf u}-{\bf u}^*_{A, h_l},\bw)|\\
  &\leq&|\bu-\bu^*_{A, h_l}|^2_{1,\Omega_f}|\bw-\bv_{h_l}|_{1,\Omega_f}
   +|\bu-\bu^*_{A, h_l}|_{1,\Omega_f}\|\bu-\bu^*_{A, h_l}\|_{0,\Omega_f}\|\bw\|_{2,\Omega_f}
\en
\bn
|c(\bu-\bu^{h_{l-1}}_A,\bu-\bu^*_{A, h_l},\bv_{h_l})|
  &\leq&|c({\bf u}-{\bu}^{h_{l-1}}_A,{\bf u}-{\bf u}^*_{A, h_l},\bw-{\bf v}_{h_l})|
    +|c({\bf u}-{\bu}^{h_{l-1}}_A,{\bf u}-{\bf u}^*_{A, h_l},\bw)|\\
  &\leq&|\bu-\bu^{h_{l-1}}_A|_{1,\Omega_f}|\bu-\bu^*_{A, h_l}|_{1,\Omega_f}|\bw-\bv_{h_l}|_{1,\Omega_f}\\
   &&+\|\bu-\bu^{h_{l-1}}_A\|_{0,\Omega_f}|\bu-\bu^*_{A, h_l}|_{1,\Omega_f}\|\bw\|_{2,\Omega_f}\\
\en
The estimate of the last term is similar to that for the above term
\bn
|c(\bu-\bu^*_{A, h_l},\bu-\bu^{h_{l-1}}_A,\bv_{h_l})|
  &\leq &  |\bu-\bu^{h_{l-1}}_A|_{1,\Omega_f}|\bu-\bu^*_{A, h_l}|_{1,\Omega_f}|\bw-\bv_{h_l}|_{1,\Omega_f}\\
   &&+\|\bu-\bu^{h_{l-1}}_A\|_{0,\Omega_f}|\bu-\bu^*_{A, h_l}|_{1,\Omega_f}\|\bw\|_{2,\Omega_f}\\
\en
Putting all the terms in the right hand side of \reff{L2-estimates-6} together, and combining with the regularity estimate, discarding the terms which are of the same order or higher order errors, we see that \reff{L2-estimates-0} holds true.
\Endproof

We comment here that Theorem \ref{M-lev-Err-A} is valid for a general $k$-th order discretization. However, because of the complex forms of the error terms, it is not easy to identify the scaling relationship for the two adjacent mesh level sizes. As previously mentioned, we are particularly interested in the first order and the second order discretizations. The following theorem states that for the first and the second order discretizations, if $h_l= h_{l-1}^2$, the energy norm errors of the final-step solution and the $L^2$ norm of $\bu^{h_l}_A$ are of the same orders as those of the FE solution.
\begin{thm}\label{M-lev-Err-Exam}
Let $(\phi,\bu,p)$, $(\phi^*_{A, h_l},\bu^*_{A,h_l},p^*_{A, h_l})$ and $(\phi^{h_l}_A,\bu^{h_l}_A,p^{h_l}_A)$ be defined by the problems \reff{Variational-Pro}, \reff{FG_darcy1_A}-\reff{FG_ns1_A}, and \reff{FG_darcy2_A}-\reff{FG_ns2_A}, respectively. For the first order and the second discretizations, i.e., $k=1$ or $k=2$, if $h_l=h_{l-1}^2$, under the assumptions of Lemma \ref{FEM-Err}, there hold
     \be\label{M-lev-Err-phi-*-Exam}
       |\phi-\phi^*_{A, h}|_{1,\Omega_p} + |\bu-{\bf u}^*_{A,h}|_{1,\Omega_f}+ \|p-p^*_{A,h}\|_{0,\Omega_f}\lesssim  h^k_l +h^{k+1}_{l-1};
     \ee
\be\label{L2-estimates-*-0-Exam}
  \|{\bf u}-{\bf u}^*_{A, h}\|_{0,\Omega_f}\lesssim h^2_l+h_l h^2_{l-1}, \quad \mbox{if}~ k=1; \quad   \|{\bf u}-{\bf u}^*_{A, h}\|_{0,\Omega_f}\lesssim h^3_l+h_l h^3_{l-1}, \quad \mbox{if}~ k=2;
\ee
\be\label{T-lev-Err-phiup-h-Exam}
       |\phi-\phi^h_A|_{1,\Omega_p} +  |\bu-{\bf u}^h_A|_{1,\Omega_f}+  \|p-p^h_A\|_{0,\Omega_f} \lesssim h^k_l;
\ee
\be\label{L2-estimates-h-0-Exam}
  \|{\bf u}-{\bf u}^h\|_{0,\Omega_f}\lesssim  h^{k+1}_l.
\ee
\end{thm}

\Proof We apply mathematical induction to the meshlevel $l$. For proving the results for the solution on meshlevel $l$, we will assume that the conclusions for the solutions on meshlevel $l-1$ hold true. From Lemma \ref{T-lev-Err}, we know that $L=1$ (in \reff{T-lev-Err-phi-*}-\reff{L2-estimates-h-0} by changing $H$ to be $h_{0}$ and $h$ to be $h_1$), the error estimates for the intermediate-step solution, and the final-step solution hold true. For both $k=1$ and $k=2$, if $h_1=h^2_0$ the estimates \reff{M-lev-Err-phi-*-Exam}-\reff{L2-estimates-h-0-Exam} hold true (see also {\bf Remark 4.1} in \cite{huang2016newton}). We are going to prove the results for a general meshlevel $l$. We will discuss the two cases: $k=1$ and $k=2$ separately.

If $k=1$, by using \reff{M-lev-Err-phi*}-\reff{L2-estimates-0},
 we see that the estimates for  $(\phi^*_{A, h_l},\bu^*_{A,h_l},p^*_{A, h_l})$ and $(\phi^{h_l}_A,\bu^{h_l}_A,p^{h_l}_A)$ are as follows.
$$
|\phi-\phi^{*}_{A, h_l}|_{1,\Omega_p}\lesssim h_{l-1}^2+h_{l-1}^2 +h_l \lesssim h_l;
$$
$$
|{\bf u}-{\bf u}^*_{A, h_l}|_{1,\Omega_f}+\|p-p^*_{A, h_l}\|_{0,\Omega_f} \lesssim (h_l +h^2_{l-1})+ h^2_{l-1} + h_l \lesssim h_l;
$$
$$
\|\bu-\bu^*_{A, h_l}\|_{0,\Omega_f}  \lesssim  h_l ( h_l + h_{l-1}^2)+( h_l + h_{l-1}^2)h^2_{l-1}+h_l^2 \lesssim h_l^2 + h_l h_{l-1}^2 ;
$$
$$
|\phi-\phi^{h_l}_A|_{1,\Omega_p} \lesssim h_l ( h_l + h_{l-1}^2) + h_l + ( h_l + h_{l-1}^2) \lesssim h_l;
$$
$$
|{\bf u}-{\bu}^{h_l}_A|_{1,\Omega_f}+\|p-p^{h_l}_A\|_{0,\Omega_f}\lesssim (h_l + h_{l-1}^2) + ( h_l + h_{l-1}^2)h_{l-1}+(h_l + h_{l-1}^2)(h_l+h_{l-1})+ h_l \lesssim h_l;
$$
$$
\|{\bf u}-{\bu}^{h_l}_A\|_{0,\Omega_f} \lesssim  h_l \big(h_l + h_{l-1}^2 + h_{l-1}( h_l + h_{l-1}^2)  \big) +h_l^2 + h^2_{l-1} (h_l + h_{l-1}^2) \lesssim h_l^2.
$$

If $k=2$, we see that we see that the estimates for  $(\phi^*_{A, h_l},\bu^*_{A,h_l},p^*_{A, h_l})$ and $(\phi^{h_l}_A,\bu^{h_l}_A,p^{h_l}_A)$ are as follows.
$$
|\phi-\phi^{*}_{A, h_l}|_{1,\Omega_p}\lesssim h_{l-1}h_{l-1}^2+h_{l-1}^3 +h^2_l \lesssim h^2_l + h_{l-1}^3;
$$
$$
|{\bf u}-{\bf u}^*_{A, h_l}|_{1,\Omega_f}+\|p-p^*_{A, h_l}\|_{0,\Omega_f} \lesssim h^2_l + h_{l-1}^3+ h^2_{l-1}h^3_{l-1} + h^2_l +h^2_l+h_{l-1}^3 \lesssim   h^2_l + h_{l-1}^3;
$$
$$
\|\bu-\bu^*_{A, h_l}\|_{0,\Omega_f}  \lesssim  h_l ( h_l^2 + h_{l-1}^3)+  (h_l^2+h_{l-1}^3) h^3_{l-1}+h_l^3 \lesssim h_l^3 + h_l h_{l-1}^3 ;
$$
$$
|\phi-\phi^{h_l}_A|_{1,\Omega_p} \lesssim  h_l ( h_l^2 + h_{l-1}^3)+ h_l^3 + h_l h^3_{l-1} + h_l^2 \lesssim h^2_l;
$$
$$
|{\bf u}-{\bu}^{h_l}_A|_{1,\Omega_f}+\|p-p^{h_l}_A\|_{0,\Omega_f}\lesssim h_l^2 + h_l h_{l-1}^3 + h_l^2 +  (h^2_l + h_{l-1}^3) h^2_{l-1}) \lesssim h^2_l;
$$
$$
\|{\bf u}-{\bu}^{h_l}_A\|_{0,\Omega_f} \lesssim  h_l \big(h^2_l + h_{l-1}^3( h^3_l + h_l h_{l-1}^3)  \big) + h^3_{l-1} (h^2_l + h_l h_{l-1}^3) +h_l^3 \lesssim h_l^3.
$$

\Endproof

From Theorem \ref{M-lev-Err-Exam}, under the scaling $h_l=h^{2}_{l-1}$, it is shown in \reff{L2-estimates-*-0-Exam} that the intermediate-step solution does not have optimal $L^2$ errors. To ensure the intermediate-step solution has optimal $L^2$ errors, one usually requires a very stringent scaling between the meshsizes of the two subsequent mesh levels. In practice, we are not interested in making the intermediate-step solution has optimal $L^2$ error. The estimate of the $L^2$ errors is for the purpose of estimating the energy norm of the final-step solution.

We would comment that the theoretical analysis of {\it Algorithm D} can be done similar to that for {\it Algorithm A}. Noting that there is no correction step in {\it Algorithm D}, the scaling of the meshsizes between two adjacent meshlevels are more stringent than that for {\it Algorithm A}. In the next section, we provide numerical experiments showing that for the first order discretization, the final-step solution of {\it Algorithm D} is still optimal if $h_l=h_{l-1}^2$. However, for the second order discretization, one has to take $h_l=h_{l-1}^{3/2}$ to ensure the final-step solution in the energy norm is optimal (in particular, for the variable $\phi^{h_l}_D$).

\section{Numerical Experiments}
We now present numerical experiments to demonstrate the effectiveness and the accuracy of the multi-level approach. In order to make our experiments more solid, we first compare different two-level algorithms then give the numerical experiments for the multilevel cases.

The computational domain is $\Omega \subset \mathbb{R}^2$ with $\Omega_f = (0, 1) \times (1, 2)$, $\Omega_p = (0, 1) \times (0, 1)$ and the interface $\Gamma= (0, 1) \times \{1\}$. The components of ${\bf u}$ are denoted by $(u, v)$. For simplicity, all the parameters in the coupled NS/Darcy model are set to $1$. The boundary conditions and right hand side functions of the coupled NS/Darcy model are chosen so that the exact solution $(u, v, p, \phi)$ is given by
\begin{equation}\label{model_problem}
\cases{ u =\mbox{cos}(\frac{\pi y}{2})^2\mbox{sin}(\frac{\pi x}{2}
), \cr v= -\mbox{cos}(\frac{\pi x}{2})(\frac{1}{4} \mbox{sin}(\pi
y)+\frac{\pi y}{4}), \cr p=\frac{\pi}{4} \mbox{cos}(\frac{\pi x}{2})
(y-1-\mbox{cos}(\pi y)), \cr {\phi} =\frac{\pi
y}{4}\mbox{cos}(\frac{1}{2} \pi x).}
\end{equation}
The coupled nonlinear FE problem is solved by the Picard iteration: given $(\ul^{0},p^{0})\in \Xl_h \times Q_h$, for $m\geq0$, find $(\ul^{m+1},p^{m+1})\in \Xl_h \times Q_h$ such that
\be\label{picard-ns-Darcy}
               \left \{
                \ba {ll}
                 a(\ul^{m+1},\vl)+c(\bu^m,\bu^{m+1},\bv)+b(\bv,p^{m+1})=f(\vl)\qquad & \forall \vl=({\bf v},\psi)\in \Xl_h,\\
                 b(\bu^{m+1},q)=0\qquad & \forall q\in Q_h.
                \ea
               \right.
\ee
The stopping criterion for the Picard iteration is $\|{\bf U}^{m+1}-{\bf U}^{m}\|_{l^2} < 10^{-7}$, where ${\bf U}^{m}$ is the nodal-value vector for the $m$-th iterate. In all algorithms, the symmetric positive definite linear system of the fine-grid Darcy problems are solved by the PCG method with the incomplete Cholesky factorization as preconditioner. The stopping criterion of PCG is set to be $10^{-9}$ and the dropping tolerance of the incomplete Cholesky factorization is $10^{-3}$. For solving the fine-grid linearized Navier-Stokes problems and the coarse-grid linear system at each step of the Picard iteration, we employ the preconditioned GMRES method with the stopping criterion $\frac{\|{\bf r}_q\|_{l^2}}{\|{\bf r}_0\|_{l^2}} < 10^{-9}$, where ${\bf r}_q$ is the residual at the $q$-th iteration of the GMRES method. The preconditioners of these saddle point problems are designed by using the Green function theory \cite{kay2002preconditioner}. Interested readers are referred to \cite{cai2008modeling, cai2012decoupled} for more details. All experiments were performed using personal desktop computer with the processor Intel Core i3 2130 (Operating speed 3.4 GHz). For the tests we presented in this paper, the average number of Picard iteration is 6, the number of GMRES iterations for the coupled model or linearized Navier-Stokes model with Green function theory based preconditioner is around 30, the iterations of PCG with incomplete Cholesky factorization preconditioning are less than 280 in all tests.

In the implementation of the two-level and multilevel algorithms, the key part is the Finite Element interpolations. FE interpolations are applied from coarse grid to fine grid or between different submodels. For example, when solving the Darcy problem (\ref{FG_darcy1_A}), we need to compute the Neumann data at the quadrature points of the fine grid by using the coarse grid solution. Standard FE interpolation is applied to supplement the Neumann data: take the coarse grid NS solution, and use the coarse grid basis functions to calculate the Neuman data at the quadrature points when assembling the right hand side of (\ref{FG_darcy1_A}).

For {\it Algorithm A}, the following notations are used to measure the solution errors, the intermediate-step two-level solution errors and the final two-level solution errors for $\phi$ in the energy norm and the $L^2$ norm.
$$
\begin{array}{ll}
e^{\phi^*_h}_{0, A}= ||{\phi}^*_{A, h}-\phi||_{0,\Omega_p}, \quad & e^{\phi^h}_{0, A}= ||{\phi}^h_A-\phi||_{0,\Omega_p},  \cr
e^{\phi^*_h}_{1, A}= |{\phi}^*_{A, h}-\phi|_{1,\Omega_p},  \quad & e^{\phi^h}_{1, A}= |{\phi}^h_A-\phi|_{1,\Omega_p}.
\end{array}
$$
Similarly, the notations, $e^{u^*_h}_{0, A}, e^{u^*_h}_{1, A}, e^{u^h}_{0, A}, e^{u^h}_{1, A}, e^{v^*_h}_{0, A}, e^{v^*_h}_{1, A}, e^{v^h}_{0, A}, e^{v^h}_{1, A}, e^{p^*_h}_{0, A}$, and $e^{p^h}_{0, A}$ are used to denote the corresponding errors for the velocity components and the pressure variable with specified norms. For the coupled nonlinear FE algorithm, we use $e^{\phi_h}_{1}, e^{u_h}_{1}, e^{v_h}_1, e^{p_h}_0$ to denote the corresponding finite element errors. Similarly, if the algorithm is changed to be {\it Algorithm B}, {\it Algorithm C} or {\it Algorithm D}, the subindex of the errors will be changed correspondingly. We keep $4$ valid digits when calculating all the errors in the following tests.

\subsection{Comparisons of the two-level algorithms}

\begin{table}[ht]
\begin{center}
{\scriptsize
\begin{tabular}{|c||c|c||c|c||c|c||c|}
\hline
$h=H^3$ &$e^{\phi_h}_0$  &$e^{{\phi}_h}_1$ &$e^{u_h}_0$ &$e^{u_h}_1$ &$e^{v_h}_0$ &$e^{v_h}_1$ &$e^{p_h}_0$  \\
$1/8$  &1.736E-3 &6.134E-2 &3.685E-3 &1.263E-1 &2.588E-3 &1.066E-1 &7.420E-2 \\
$1/27$ &1.552E-4 &1.823E-2 &3.213E-4 &3.714E-2 &2.251E-4 &3.070E-2 &9.113E-3 \\
$1/64$ &2.766E-5 &7.693E-3 &5.697E-5 &1.564E-2 &3.996E-5 &1.289E-2 &2.255E-3 \\
$1/125$&7.253E-6 &3.939E-3 &1.491E-5 &8.000E-3 &1.046E-5 &6.587E-3 &7.906E-4  \\
\hline
\hline
$h=H^3$ &$e^{\phi^*_h}_{0, A}$  &$e^{\phi^*_h}_{1, A}$ &$e^{u^*_h}_{0, A}$ &$e^{u^*_h}_{1, A}$ &$e^{v^*_h}_{0, A}$ &$e^{v^*_h}_{1, A}$ &$e^{p^*_h}_{0, A}$  \\
$1/8$  &7.649E-3 &6.848E-2 &3.830E-3 &1.263E-1 &2.573E-3 &1.067E-1 &7.541E-2 \\
$1/27$ &3.051E-3 &2.276E-2 &4.089E-4 &3.718E-2 &2.392E-4 &3.072E-2 &1.355E-2 \\
$1/64$ &1.646E-3 &1.073E-2 &9.572E-5 &1.567E-2 &6.792E-5 &1.290E-2 &5.606E-3 \\
$1/125$&1.049E-3 &6.206E-3 &4.773E-5 &8.024E-3 &3.034E-5 &6.600E-3 &3.026E-3 \\
$h=H^3$ &$e^{\phi^h}_{0, A}$  &$e^{\phi^h}_{1, A}$ &$e^{u^h}_{0, A}$ &$e^{u^h}_{1, A}$ &$e^{v^h}_{0, A}$ &$e^{v^h}_{1, A}$ &$e^{p^h}_{0, A}$  \\
$1/8$  &1.741E-3 &6.134E-2 &3.685E-3 &1.263E-1 &2.588E-3 &1.066E-1 &7.421E-2 \\   
$1/27$ &1.580E-4 &1.823E-2 &2.891E-4 &3.715E-2 &2.251E-4 &3.070E-2 &9.194E-3 \\
$1/64$ &2.766E-5 &7.693E-3 &5.156E-5 &1.564E-2 &3.988E-5 &1.289E-2 &2.262E-3 \\
$1/125$&6.686E-6 &3.939E-3 &1.361E-5 &8.000E-3 &1.044E-5 &6.587E-3 &7.919E-4 \\
\hline
\hline
$h=H^3$ &$e^{\phi^*_h}_{0, B}$  &$e^{\phi^*_h}_{1, B}$ &$e^{u^*_h}_{0, B}$ &$e^{u^*_h}_{1, B}$ &$e^{v^*_h}_{0, B}$ &$e^{v^*_h}_{1, B}$ &$e^{p^*_h}_{0, B}$  \\
$1/8$  &1.725E-3 &6.134E-2 &3.692E-3 &1.263E-1 &2.573E-3 &1.067E-1 &9.835E-2 \\
$1/27$ &1.476E-4 &1.823E-2 &4.163E-4 &3.716E-2 &2.633E-4 &3.075E-2 &3.286E-2 \\
$1/64$ &2.195E-5 &7.693E-3 &1.206E-4 &1.567E-2 &1.143E-4 &1.293E-2 &1.766E-2 \\
$1/125$&6.615E-6 &3.939E-3 &7.639E-5 &8.031E-3 &8.135E-5 &6.634E-3 &1.111E-3 \\
$h=H^3$ &$e^{\phi^h}_{0, B}$  &$e^{\phi^h}_{1, B}$ &$e^{u^h}_{0, B}$ &$e^{u^h}_{1, B}$ &$e^{v^h}_{0, B}$ &$e^{v^h}_{1, B}$ &$e^{p^h}_{0, B}$  \\
$1/8$  &1.736E-3 &6.134E-2 &3.685E-3 &1.263E-1 &2.588E-3 &1.066E-1 &7.418E-2 \\
$1/27$ &1.553E-4 &1.823E-2 &3.225E-4 &3.714E-2 &2.252E-4 &3.070E-2 &9.110E-3 \\
$1/64$ &2.773E-5 &7.693E-3 &5.767E-5 &1.564E-2 &3.998E-5 &1.289E-2 &2.248E-3 \\
$1/125$&7.319E-6 &3.939E-3 &1.363E-5 &8.000E-3 &1.052E-5 &6.587E-3 &7.879E-4 \\
\hline
\hline
$h=H^3$ &$e^{\phi^h}_{0, C}$  &$e^{\phi^h}_{1, C}$ &$e^{u^h}_{0, C}$ &$e^{u^h}_{1, C}$ &$e^{v^h}_{0, C}$ &$e^{v^h}_{1, C}$ &$e^{p^h}_{0, B}$  \\
$1/8$  &7.649E-3 &6.848E-2 &3.692E-3 &1.263E-1 &2.573E-3 &1.067E-1 &9.835E-2 \\
$1/27$ &3.051E-3 &2.276E-2 &4.163E-4 &3.716E-2 &2.633E-4 &3.075E-2 &3.286E-2 \\
$1/64$ &1.646E-3 &1.073E-2 &1.206E-4 &1.567E-2 &1.143E-4 &1.293E-2 &1.766E-2 \\
$1/125$&1.049E-3 &6.206E-3 &7.639E-5 &8.031E-3 &8.135E-5 &6.634E-3 &1.111E-3 \\
\hline
\end{tabular}
}
\caption{Tests for the two-level algorithms under the Mini/$P_1$ element discretizations. The FE solution errors, the intermediate-step two-level solution errors, and the final two-level solution errors with the scaling between the two-level sizes being set as $h=H^3$.
}
\end{center}
\label{Mini-TG-Err-H3}
\end{table}

\begin{figure}[ht]
\centering
\includegraphics[height=5.6cm]{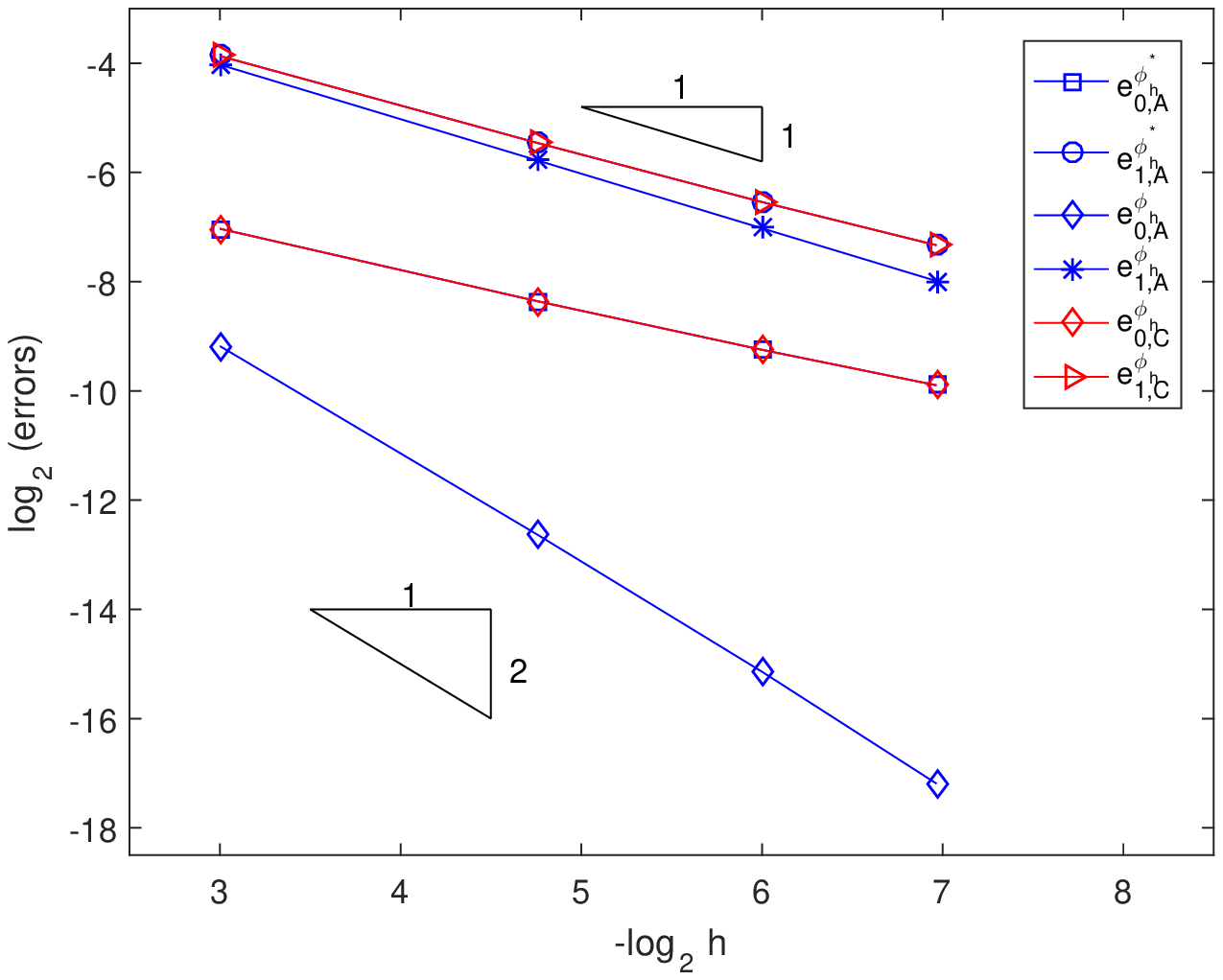}
\hskip -0.2cm
\includegraphics[height=5.6cm]{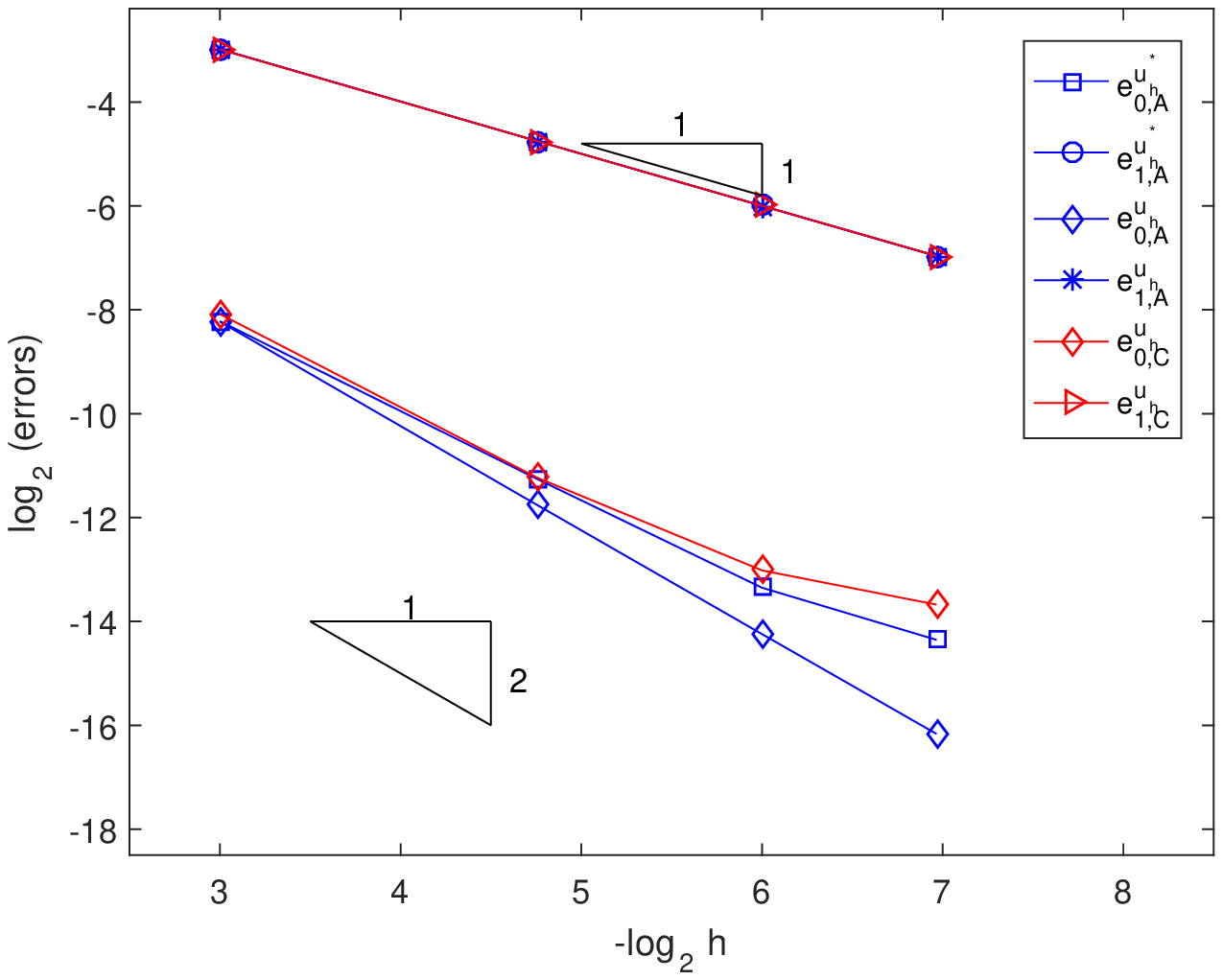}
\vskip .1cm
\includegraphics[height=5.6cm]{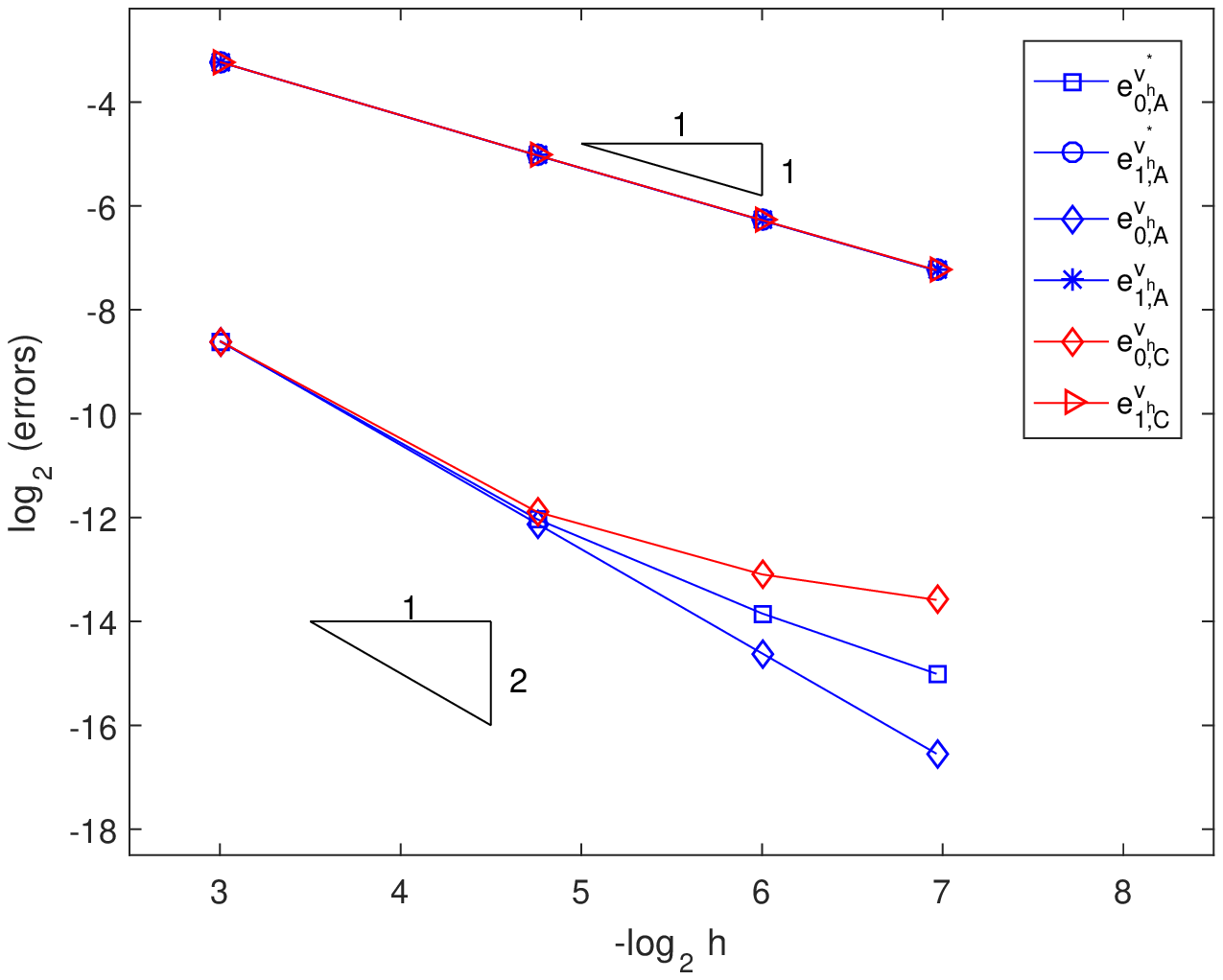}
\hskip -0.2cm
\includegraphics[height=5.6cm]{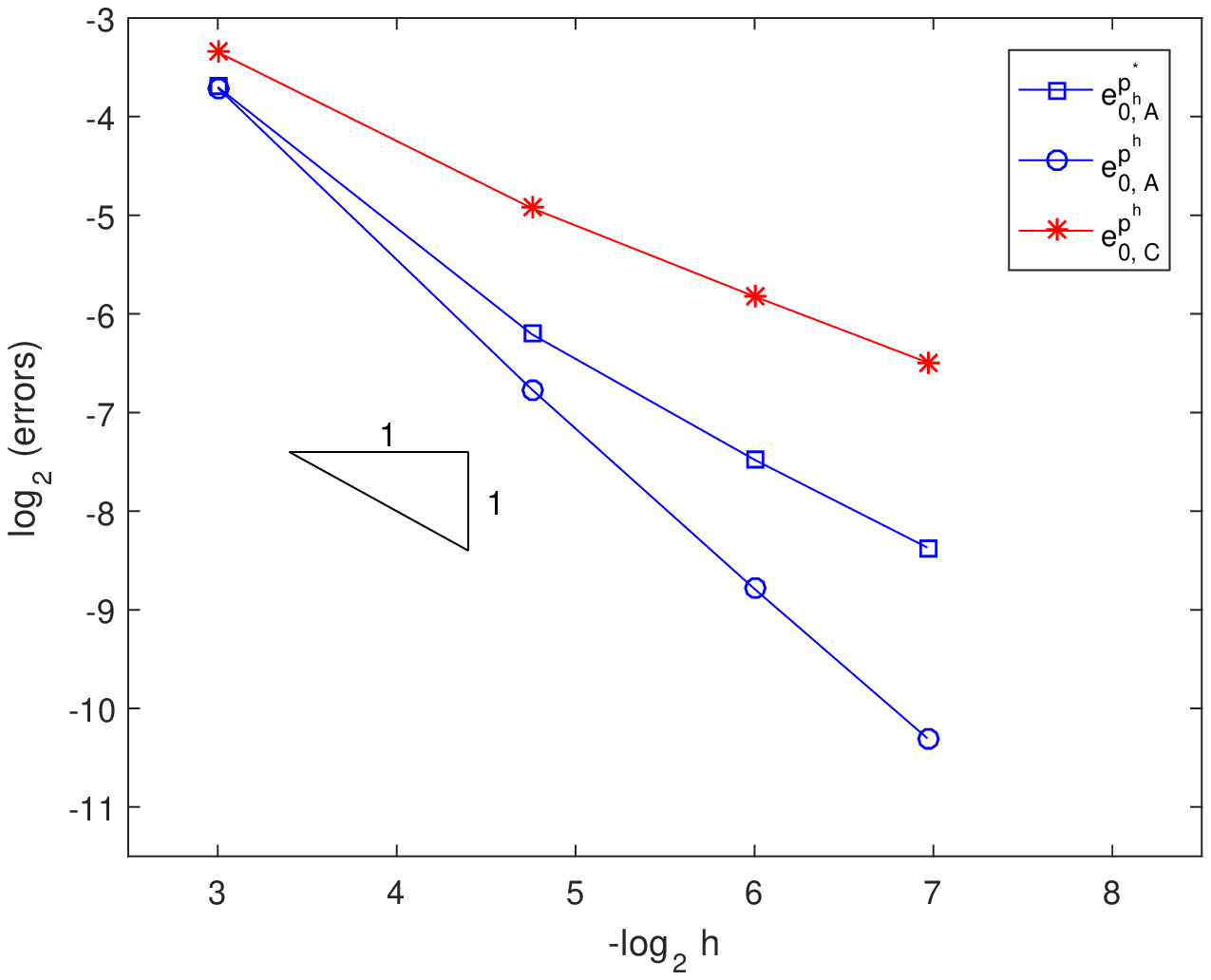}
\vspace{-.2cm}
\caption{Plots of the errors from Table \ref{Mini-TG-Err-H3}. Comparisons of {\it Algorithm A}
 and {\it Algorithm C}. 
}
\label{fig-Mini-TL-error}
\end{figure}

\begin{table}[ht]
\begin{center}
{\scriptsize
\begin{tabular}{|r||c|c||c|c||c|c||c|}
\hline
$h \approx H^{5/2}, H$ &$e^{\phi_h}_0$  &$e^{{\phi}_h}_1$ &$e^{u_h}_0$ &$e^{u_h}_1$ &$e^{v_h}_0$ &$e^{v_h}_1$ &$e^{p_h}_0$  \\
$1/6, 1/2$  &1.056E-4 &2.522E-3 &4.017E-4 &8.270E-3 &2.266E-4 &4.882E-3 &2.837E-3 \\
$1/16, 1/3$ &5.584E-6 &3.648E-4 &2.221E-5 &1.156E-3 &1.161E-5 &6.732E-4 &2.930E-4  \\
$1/32, 1/4$ &7.000E-7 &9.201E-5 &2.797E-6 &2.888E-4 &1.446E-6 &1.678E-4 &7.002E-5 \\
$1/56, 1/5$ &1.308E-7 &3.016E-5 &5.228E-7 &9.432E-5 &2.696E-7 &5.477E-5 &2.263E-5  \\
\hline
\hline
$h \approx H^{5/2}, H$ &$e^{\phi^*_h}_{0, A}$  &$e^{\phi^*_h}_{1, A}$ &$e^{u^*_h}_{0, A}$ &$e^{u^*_h}_{1, A}$ &$e^{v^*_h}_{0, A}$ &$e^{v^*_h}_{1, A}$ &$e^{p^*_h}_{0, A}$  \\
$1/6, 1/2$  &1.808E-4 &2.830E-3 &4.006E-4 &8.270E-3 &2.275E-4 &4.883E-3 &2.849E-3 \\
$1/16, 1/3$ &3.077E-5 &4.666E-4 &2.233E-5 &1.156E-3 &1.212E-5 &6.736E-4 &2.947E-4 \\
$1/32, 1/4$ &8.907E-6 &1.317E-4 &2.955E-6 &2.889E-4 &1.794E-6 &1.680E-4 &7.060E-5 \\
$1/56, 1/5$ &3.484E-6 &5.018E-5 &6.489E-7 &9.433E-5 &4.811E-7 &5.485E-5 &2.300E-5 \\
$h \approx H^{5/2}, H$ &$e^{\phi^h}_{0, A}$  &$e^{\phi^h}_{1, A}$ &$e^{u^h}_{0, A}$ &$e^{u^h}_{1, A}$ &$e^{v^h}_{0, A}$ &$e^{v^h}_{1, A}$  &$e^{p^h}_{0, A}$  \\
$1/6, 1/2$  &1.055E-4 &2.522E-3 &4.018E-4 &8.270E-3 &2.266E-4 &4.882E-3 &2.837E-3 \\
$1/16, 1/3$ &5.587E-5 &3.648E-4 &2.221E-5 &1.156E-3 &1.161E-5 &6.732E-4 &2.930E-4 \\
$1/32, 1/4$ &7.087E-7 &9.202E-5 &2.797E-6 &2.888E-4 &1.446E-6 &1.678E-4 &7.002E-5 \\
$1/56, 1/5$ &1.380E-7 &3.016E-5 &5.229E-7 &9.432E-5 &2.696E-7 &5.477E-5 &2.263E-5 \\
\hline
\hline
$h \approx H^{5/2}, H$ &$e^{\phi^*_h}_{0, B}$  &$e^{\phi^*_h}_{1, B}$ &$e^{u^*_h}_{0, B}$ &$e^{u^*_h}_{1, B}$ &$e^{v^*_h}_{0, B}$ &$e^{v^*_h}_{1, B}$ &$e^{p^*_h}_{0, B}$  \\
$1/6, 1/2$  &1.052E-4 &2.523E-3 &4.038E-4 &8.277E-3 &2.530E-4 &4.949E-3 &3.018E-3 \\
$1/16, 1/3$ &6.126E-6 &3.656E-4 &3.195E-5 &1.160E-3 &2.663E-5 &6.978E-4 &3.472E-4 \\
$1/32, 1/4$ &1.020E-6 &9.227E-5 &7.704E-6 &2.907E-4 &7.129E-6 &1.780E-4 &9.291E-5 \\
$1/56, 1/5$ &3.095E-7 &3.027E-5 &2.916E-6 &9.535E-5 &2.697E-6 &6.035E-5 &3.489E-5 \\
$h \approx H^{5/2}, H$ &$e^{\phi^h}_{0, B}$  &$e^{\phi^h}_{1, B}$ &$e^{u^h}_{0, B}$ &$e^{u^h}_{1, B}$ &$e^{v^h}_{0, B}$ &$e^{v^h}_{1, B}$ &$e^{p^h}_{0, B}$  \\
$1/6, 1/2$  &1.057E-4 &2.522E-3 &4.019E-4 &8.270E-3 &2.266E-4 &4.882E-3 &2.836E-3 \\
$1/16, 1/3$ &5.585E-6 &3.648E-4 &2.223E-5 &1.156E-3 &1.162E-5 &6.732E-4 &2.930E-4 \\
$1/32, 1/4$ &7.002E-7 &9.201E-5 &2.800E-6 &2.888E-4 &1.449E-6 &1.678E-4 &7.003E-5 \\
$1/56, 1/5$ &1.309E-7 &3.016E-5 &5.242E-7 &9.432E-5 &2.717E-7 &5.477E-5 &2.264E-5 \\
\hline
\hline
$h \approx H^{5/2}, H$ &$e^{\phi^h}_{0, C}$  &$e^{\phi^h}_{1, C}$ &$e^{u^h}_{0, C}$ &$e^{u^h}_{1, C}$ &$e^{v^h}_{0, C}$ &$e^{v^h}_{1, C}$ &$e^{p^h}_{0, C}$  \\
$1/6, 1/2$  &1.808E-4 &2.830E-3 &4.038E-3 &8.277E-3 &2.530E-4 &4.949E-3 &3.018E-3 \\
$1/16, 1/3$ &3.077E-5 &4.666E-4 &3.195E-5 &1.160E-3 &2.663E-5 &6.978E-4 &3.472E-4 \\
$1/32, 1/4$ &8.907E-6 &1.317E-4 &7.704E-6 &2.907E-4 &7.129E-6 &1.780E-4 &9.291E-5 \\
$1/56, 1/5$ &3.484E-6 &5.018E-5 &2.916E-6 &9.535E-5 &2.697E-6 &6.035E-5 &3.489E-5 \\
\hline
\end{tabular}
}
\caption{Tests for the two-level algorithms under the Taylor-Hood/$P_2$ element discretizations. The FE solution errors, the intermediate-step two-level solution errors, and the final two-level solution errors with the scaling between the two-level mesh sizes being set as $h \approx H^{5/2}$.
}
\end{center}
\label{THP2-TG-Err-H52}
\end{table}

\begin{figure}[ht]
\centering
\includegraphics[height=5.6cm]{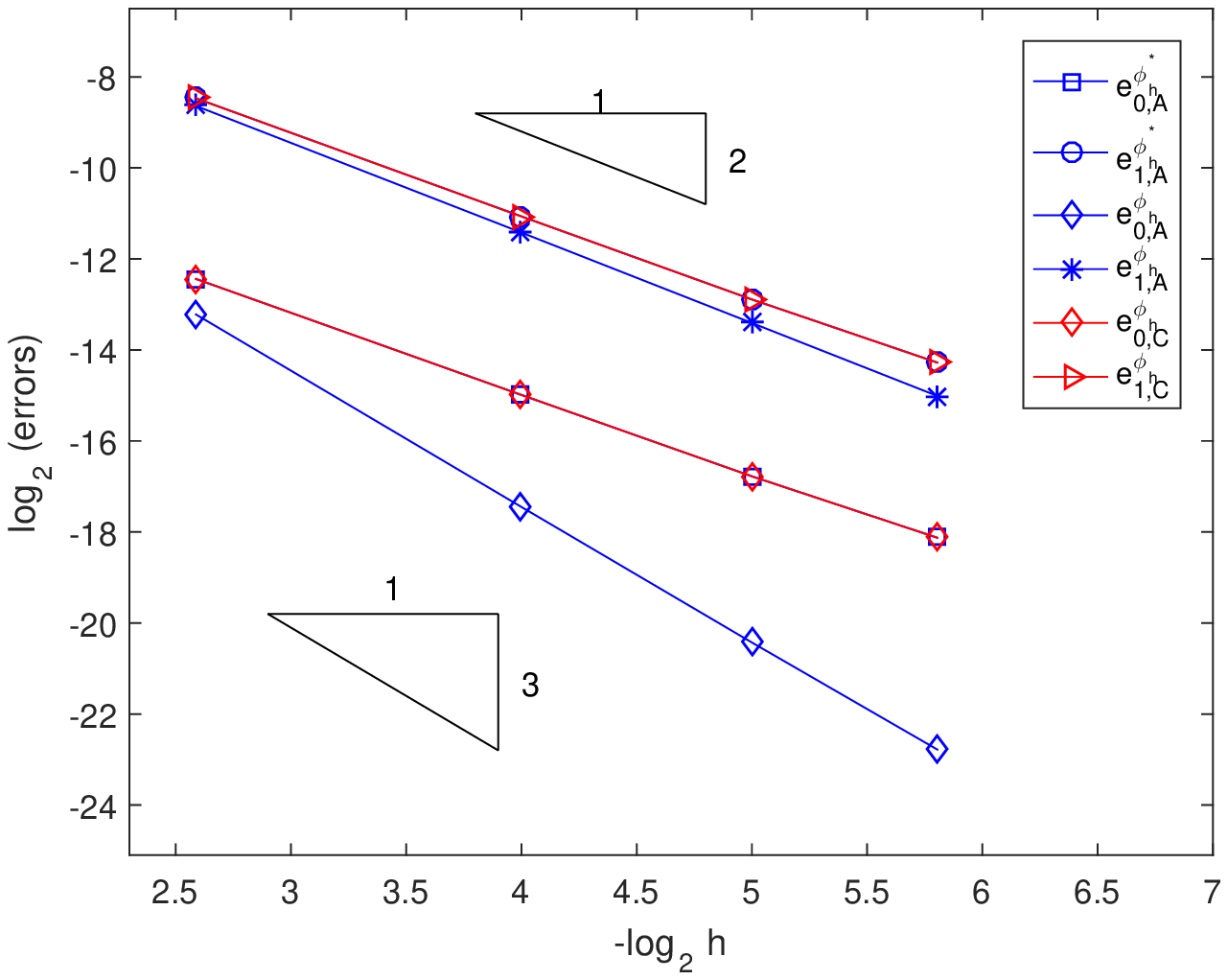}
\hskip -0.2cm
\includegraphics[height=5.6cm]{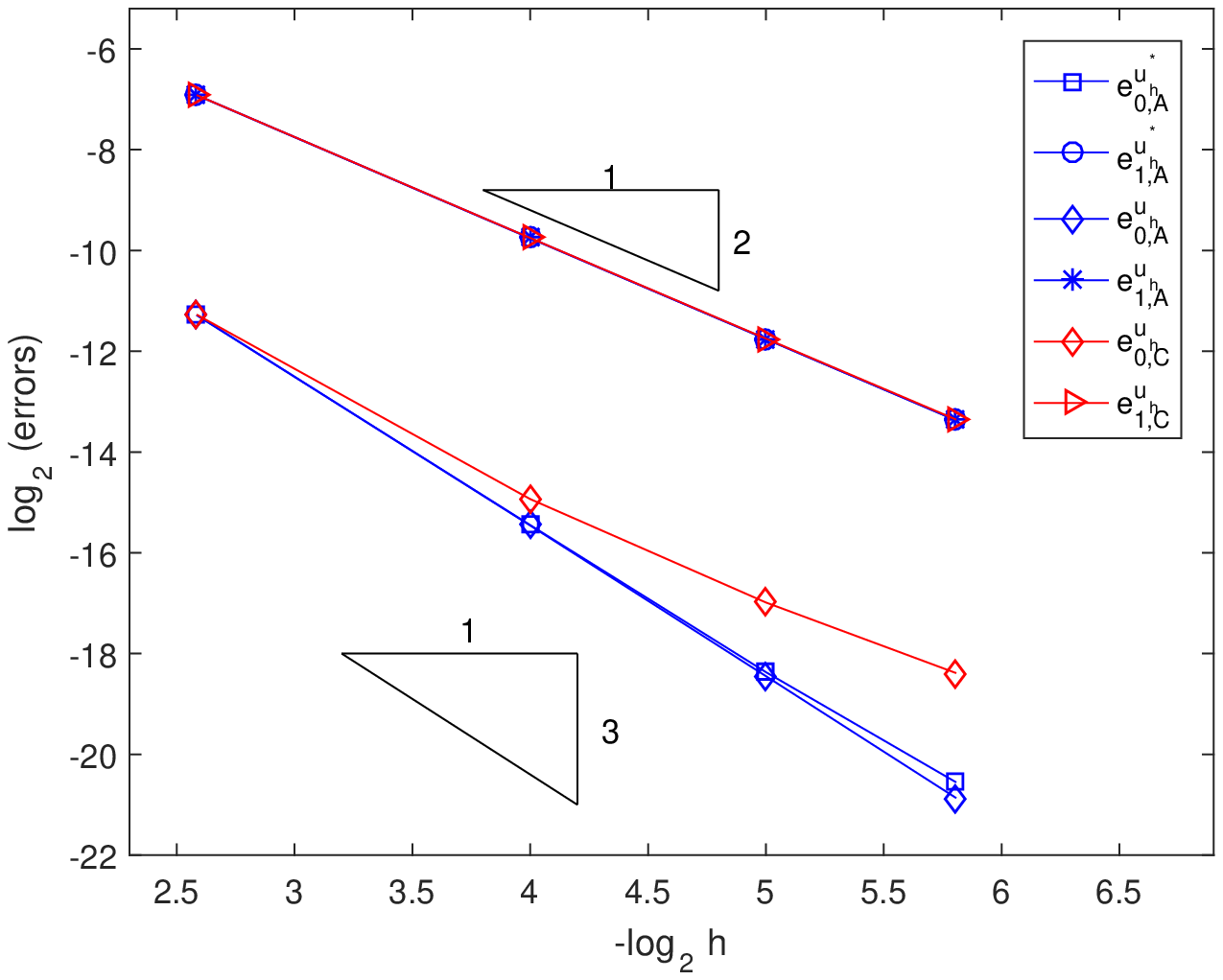}
\vskip -.05cm
\includegraphics[height=5.6cm]{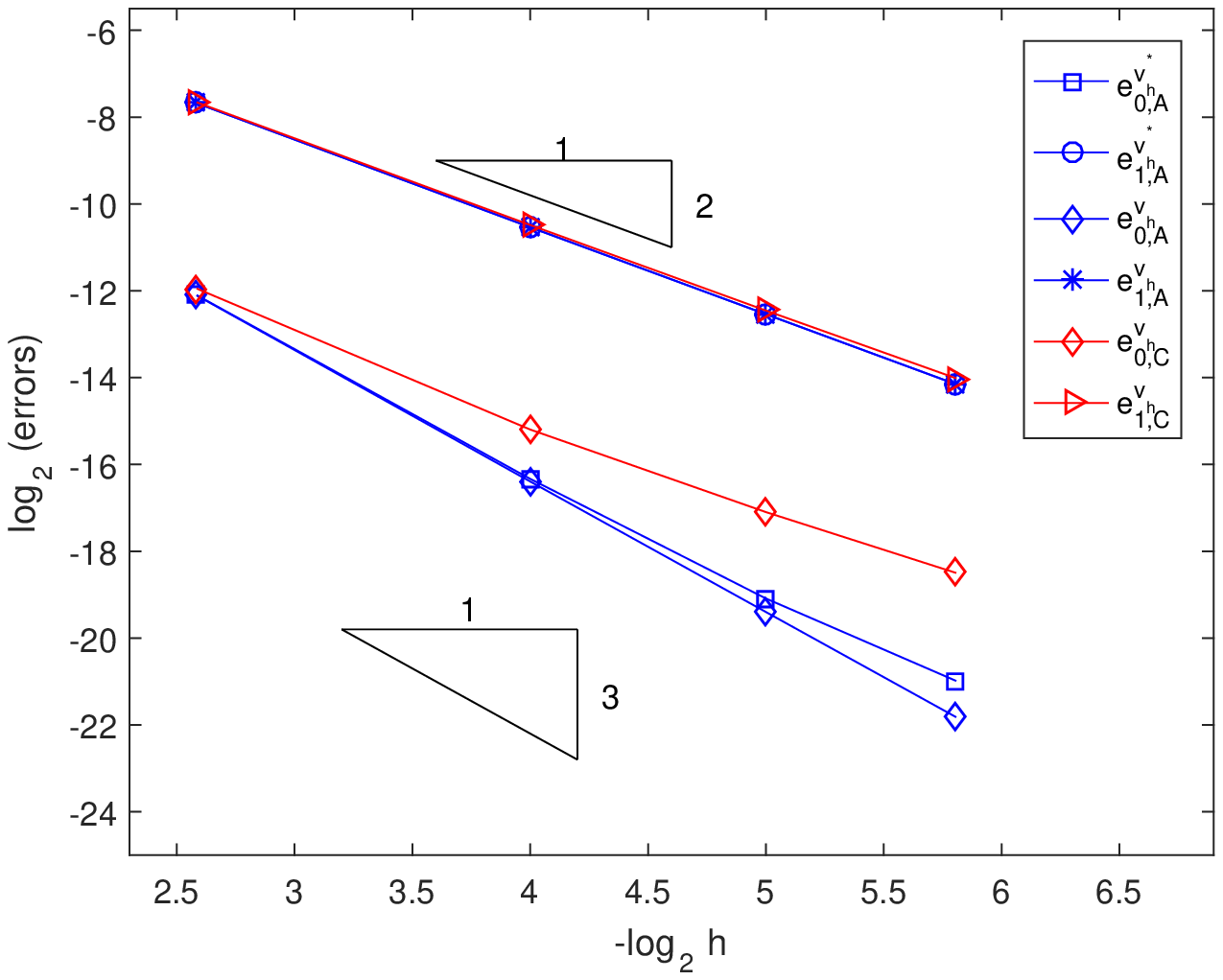}
\hskip -0.2cm
\includegraphics[height=5.6cm]{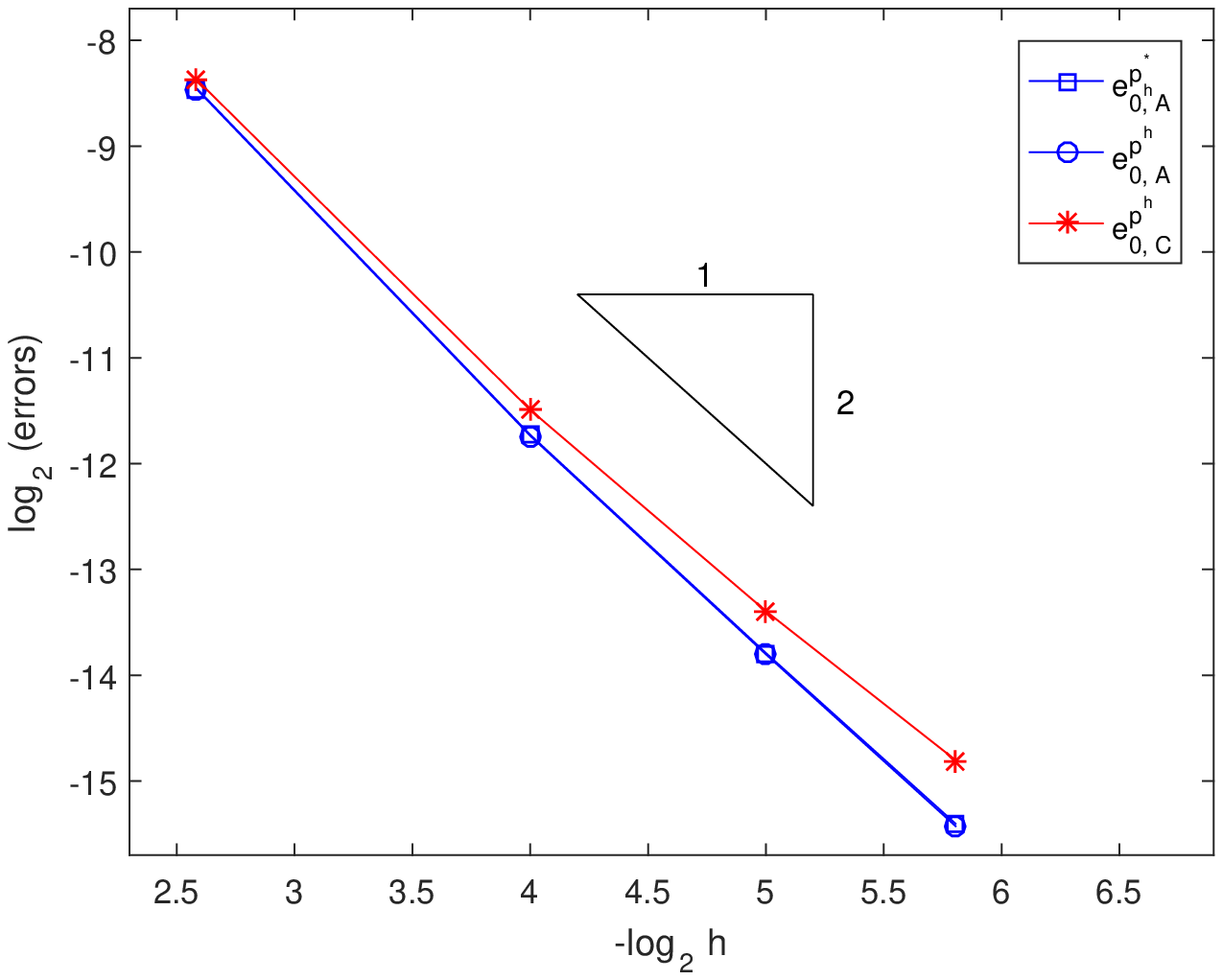}
\vspace{-.2cm}
\caption{Plots of the errors from Table 5.2. 
Comparisons of {\it Algorithm A} and {\it Algorithm C}. 
}
\label{Fig-TH-2L-error}
\end{figure}

We firstly conduct numerical experiments for comparing all algorithms under the two-level cases. According to Lemma \ref{T-lev-Err}, for the first (second) order discretization, if $h = H^3$ ($h=H^{5/2}$), then {\it Algorithm A} still gives optimal errors in the energy norm. However, it is not clear whether {\it Algorithm B} and {\it Algorithm C} also give optimal errors under such a scaling. It is also important to know the differences of the different algorithms so that we can have better understanding of their generalizations to the multilevel cases.

For the first order discretization, Mini elements and piecewise linear elements are used in $\Omega_f$ and $\Omega_p$, respectively.  The scaling between $h$ and $H$ is set to be $h=H^3$. In Table \ref{Mini-TG-Err-H3}, we report the numerical results for {\it Algorithm A}, {\it Algorithm B} and {\it Algorithm C}. For {\it Algorithm C}, we note that the results for $\phi$ actually are the same as those for $\phi$ in the intermediate-step of {\it Algorithm A} and the results for fluid variables actually are the same as those for fluid variable in the intermediate-step solution of {\it Algorithm B}. As observed from Table \ref{Mini-TG-Err-H3}, the Finite Element solution errors confirm the theoretical predictions (cf. Lemma \ref{FEM-Err}). The pressure FE error is of an order between $O(h)$ and $O(h^2)$ as there is a bubble function in the Mini element discretization for the fluid part \cite{cai2009numerical, cai2008modeling, chidyagwai2011two}. From Table \ref{Mini-TG-Err-H3}, for both {\it Algorithm A} and {\it Algorithm B}, we see that the final two-level solution errors in the energy norms, i.e., $e^{\phi^h}_{1, A}$, $e^{u^h}_{1, A}$, $e^{v^h}_{1, A}$, $e^{p^h}_{0, A}$, $e^{\phi^h}_{1, B}$, $e^{u^h}_{1, B}$, $e^{v^h}_{1, B}$, and $e^{p^h}_{0, B}$ are comparable with those of the coupled algorithm with the same meshsizes; Moreover, the $L^2$ errors of $\phi$, $u$ and $v$ are almost of the same order as that for $e^{\phi_h}_0$, $e^{u_h}_0$ and $e^{v_h}_0$. This indicates that the computational $L^2$ errors for the velocity components seem to be better than the theoretical predictions (cf. {\bf Remark 4.1} in \cite{huang2016newton}). In comparison, the intermediate-step errors $e^{\phi^*_h}_{0, A}, e^{\phi^*_h}_{1, A}, e^{u^*_h}_{0, A}$ and $e^{v^*_h}_{0, A}$ (and also $e^{\phi^*_h}_{0, B}, e^{\phi^*_h}_{1, B}, e^{u^*_h}_{0, B}$ and $e^{v^*_h}_{0, B}$) for the intermediate-step two-level solution are not optimal; the errors $e^{u^*_h}_{1, A}, e^{v^*_h}_{1, A}$ and $e^{p^*_h}_{0, A}$ (and also $e^{u^*_h}_{1, B}, e^{v^*_h}_{1, B}$ and $e^{p^*_h}_{0, B}$) are slightly worse than those corresponding errors for $u_h$ or $u^h$, $v_h$ or $v^h$ and $p_h$ or $p^h$. From the digital comparisons in Table \ref{Mini-TG-Err-H3}, we see that {\it Algorithm B} gives almost the same errors as {\it Algorithm A} in the final step. However, from Table 5.1 and Figure \ref{fig-Mini-TL-error}, {\it Algorithm C} does not give optimal error order under the same scaling setting for the two level meshsizes (in particular, for the pressure errors).


For second order discretization, Taylor-Hood elements are applied in $\Omega_f$ and piecewise quadratic elements are applied in $\Omega_p$. The scaling between $h$ and $H$ is set to be $h\approx H^{5/2}$. (Actually, except the case $h=1/32=(1/4)^{5/2}$, the fine grid sizes $h$ are selected even slightly smaller than $H^{5/2}$.) Numerical results are reported in Table 5.2 
and the comparisons of {\it Algorithm A} and {\it Algorithm C} are presented in Figure \ref{Fig-TH-2L-error}. As observed from Table 5.2 
and Figure \ref{Fig-TH-2L-error}, the Finite Element solution errors confirm the theoretical analysis of Lemma \ref{FEM-Err}; The final-step solution errors of the two-level {\it Algorithm A} in the energy norms are almost the same as those of the coupled nonlinear FE algorithm with the same meshsizes; Again, {\it Algorithm A} and {\it Algorithm B} give almost the same numerical solution; The results of {\it Algorithm C} actually correspond to the intermediate-step solution of {\it Algorithm A} and {\it Algorithm B}. From Table 5.2 and Figure \ref{Fig-TH-2L-error}, {\it Algorithm C} does not give optimal numerical errors under the scaling $h\approx H^{5/2}$ (in particular, for pressure errors). By comparing the digital results of {\it Algorithm A} and {\it Algorithm C} in Table 5.2, one can get the same conclusions as we have drawn for the first order discretization.

\subsection{Comparisons of the multilevel algorithms}

\begin{table}[ht]
\begin{center}
{\scriptsize
\begin{tabular}{|c||c|c||c|c||c|c||c|}
\hline
$h_l=h_{l-1}^2$   &$e^{\phi_h}_0$  &$e^{{\phi}_h}_1$ &$e^{u_h}_0$ &$e^{u_h}_1$ &$e^{v_h}_0$ &$e^{v_h}_1$ &$e^{p_h}_0$  \\
$2^{-1}$   &2.153E-2 &2.351E-1 &5.524E-2 &5.087E-1 &4.295E-2 &5.252E-1 &1.024E-0 \\
$2^{-2}$   &6.566E-3 &1.215E-1 &1.467E-2 &2.539E-1 &1.049E-2 &2.261E-1 &2.649E-1 \\
$2^{-4}$   &4.406E-4 &3.075E-2 &9.178E-4 &6.282E-2 &6.429E-4 &5.217E-3 &2.201E-2 \\
${2^{-8}}$ &1.729E-6 &1.923E-3 &3.551E-6 &3.905E-3 &2.493E-6 &3.214E-3 &2.628E-4 \\
\hline
\hline
$h_l=h_{l-1}^2$  &$e^{\phi^*_h}_{0, A}$  &$e^{\phi^*_h}_{1, A}$ &$e^{u^*_h}_{0, A}$ &$e^{u^*_h}_{1, A}$ &$e^{v^*_h}_{0, A}$ &$e^{v^*_h}_{1, A}$ &$e^{p^*_h}_{0, A}$  \\
${2^{-2}}$ &1.049E-2 &1.247E-1 &1.474E-2 &2.538E-1 &1.047E-2 &2.261E-1 &2.612E-1 \\
$2^{-4}$   &1.868E-3 &3.163E-2 &9.574E-4 &6.283E-2 &6.481E-4 &5.217E-3 &2.213E-2 \\
${2^{-8}}$ &1.031E-4 &1.979E-3 &3.112E-6 &3.905E-3 &3.299E-6 &3.214E-3 &3.692E-4 \\
$h_l=h_{l-1}^2$  &$e^{\phi^h}_{0, A}$  &$e^{\phi^h}_{1, A}$ &$e^{u^h}_{0, A}$ &$e^{u^h}_{1, A}$ &$e^{v^h}_{0, A}$ &$e^{v^h}_{1, A}$ &$e^{p^h}_{0, A}$  \\
${2^{-2}}$ &6.570E-3 &1.215E-1 &1.467E-2 &2.539E-1 &1.049E-2 &2.261E-1 &2.649E-1 \\
$2^{-4}$   &4.424E-4 &3.075E-2 &9.184E-4 &6.282E-2 &6.429E-4 &5.217E-3 &2.201E-2 \\
${2^{-8}}$ &1.657E-6 &1.923E-3 &3.818E-6 &3.905E-3 &2.491E-6 &3.214E-3 &2.692E-4 \\
\hline
\hline
$h_l=h_{l-1}^2$  &$e^{\phi^*_h}_{0, D}$  &$e^{\phi^*_h}_{1, D}$ &$e^{u^*_h}_{0, D}$ &$e^{u^*_h}_{1, D}$ &$e^{v^*_h}_{0, D}$ &$e^{v^*_h}_{1, D}$ &$e^{p^*_h}_{0, D}$  \\
${2^{-2}}$ &1.049E-2 &1.247E-1 &1.474E-2 &2.538E-1 &1.047E-2 &2.261E-1 &2.612E-1 \\
$2^{-4}$   &1.865E-3 &3.163E-2 &9.571E-4 &6.283E-2 &6.481E-4 &5.217E-3 &2.212E-2 \\
${2^{-8}}$ &9.986E-5 &1.976E-3 &3.177E-6 &3.905E-3 &3.332E-6 &3.214E-3 &3.507E-4 \\
$h_l=h_{l-1}^2$  &$e^{\phi^h}_{0, D}$  &$e^{\phi^h}_{1, D}$ &$e^{u^h}_{0, D}$ &$e^{u^h}_{1, D}$ &$e^{v^h}_{0, D}$ &$e^{v^h}_{1, D}$ &$e^{p^h}_{0, D}$  \\
${2^{-2}}$ &         &         &1.466E-2 &2.539E-1 &1.049E-2 &2.261E-1 &2.706E-1 \\
$2^{-4}$   &         &         &9.101E-4 &6.282E-2 &6.404E-4 &5.217E-3 &2.357E-2 \\
${2^{-8}}$ &         &         &3.318E-6 &3.905E-3 &2.746E-6 &3.214E-3 &3.933E-4 \\
\hline
\hline
$h_l=h_{l-1}^2$ &$e^{\phi^h}_{0, C}$  &$e^{\phi^h}_{1, C}$ &$e^{u^h}_{0, C}$ &$e^{u^h}_{1, C}$ &$e^{v^h}_{0, C}$ &$e^{v^h}_{1, C}$ &$e^{p^h}_{0, C}$  \\
${2^{-2}}$ &1.049E-2 &1.247E-1 &1.466E-2 &2.539E-1 &1.049E-2 &2.261E-1 &2.752E-1 \\
$2^{-4}$   &1.868E-3 &3.163E-2 &8.931E-4 &6.284E-2 &6.520E-4 &5.220E-3 &3.974E-2 \\
${2^{-8}}$ &8.642E-5 &1.967E-3 &2.998E-5 &3.910E-3 &3.463E-6 &3.233E-3 &5.265E-4 \\
\hline
\end{tabular}
}
\caption{Tests for the multilevel algorithms under the  Mini/$P_1$ element discretization. Comparisons of {\it Algorithm A}, {\it Algorithm C} and {\it Algorithm D}.
}
\end{center}
\label{Mini-ML-Err}
\end{table}

\begin{figure}[ht]
\centering
\includegraphics[height=5.6cm]{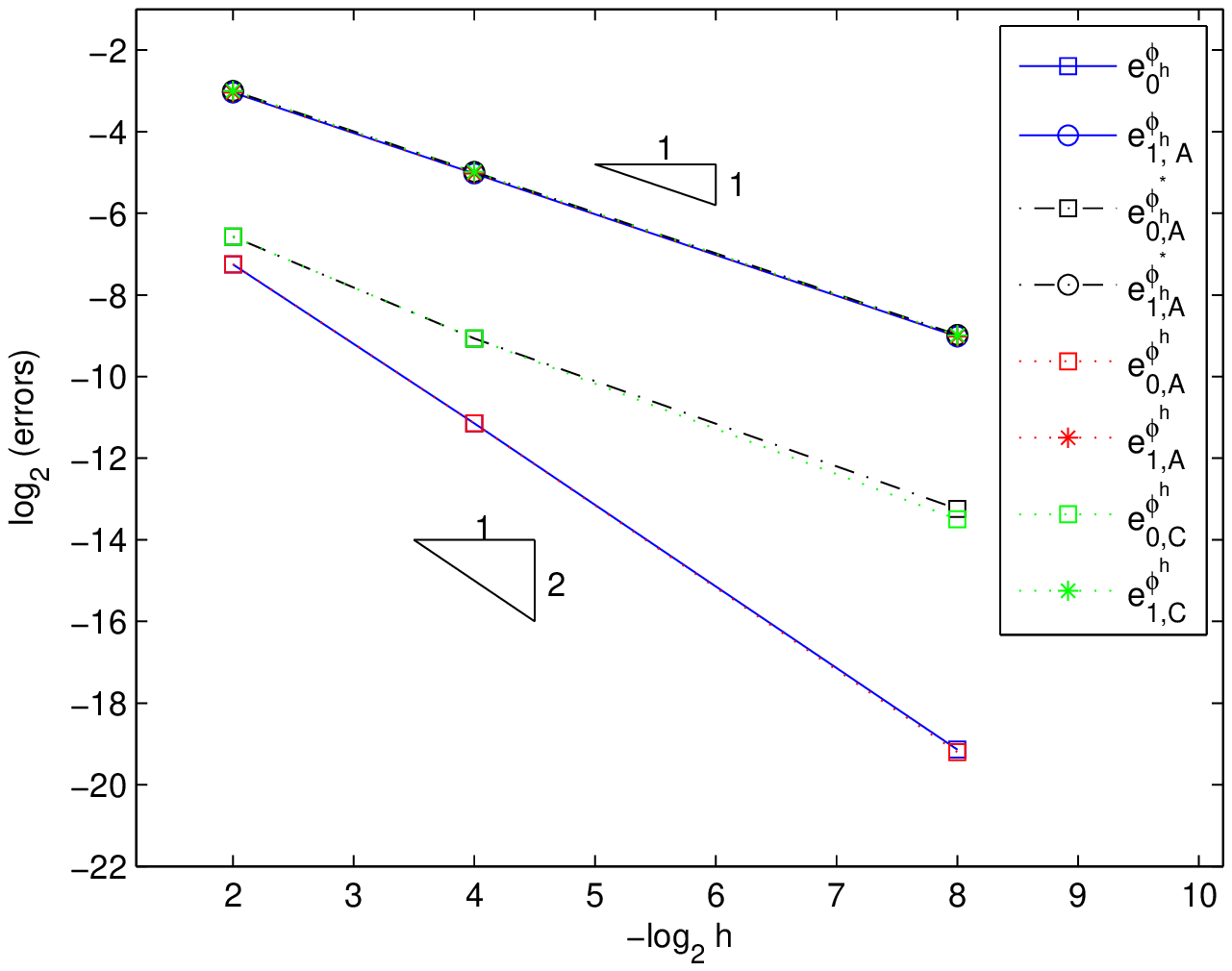}
\hskip -0.2cm
\includegraphics[height=5.6cm]{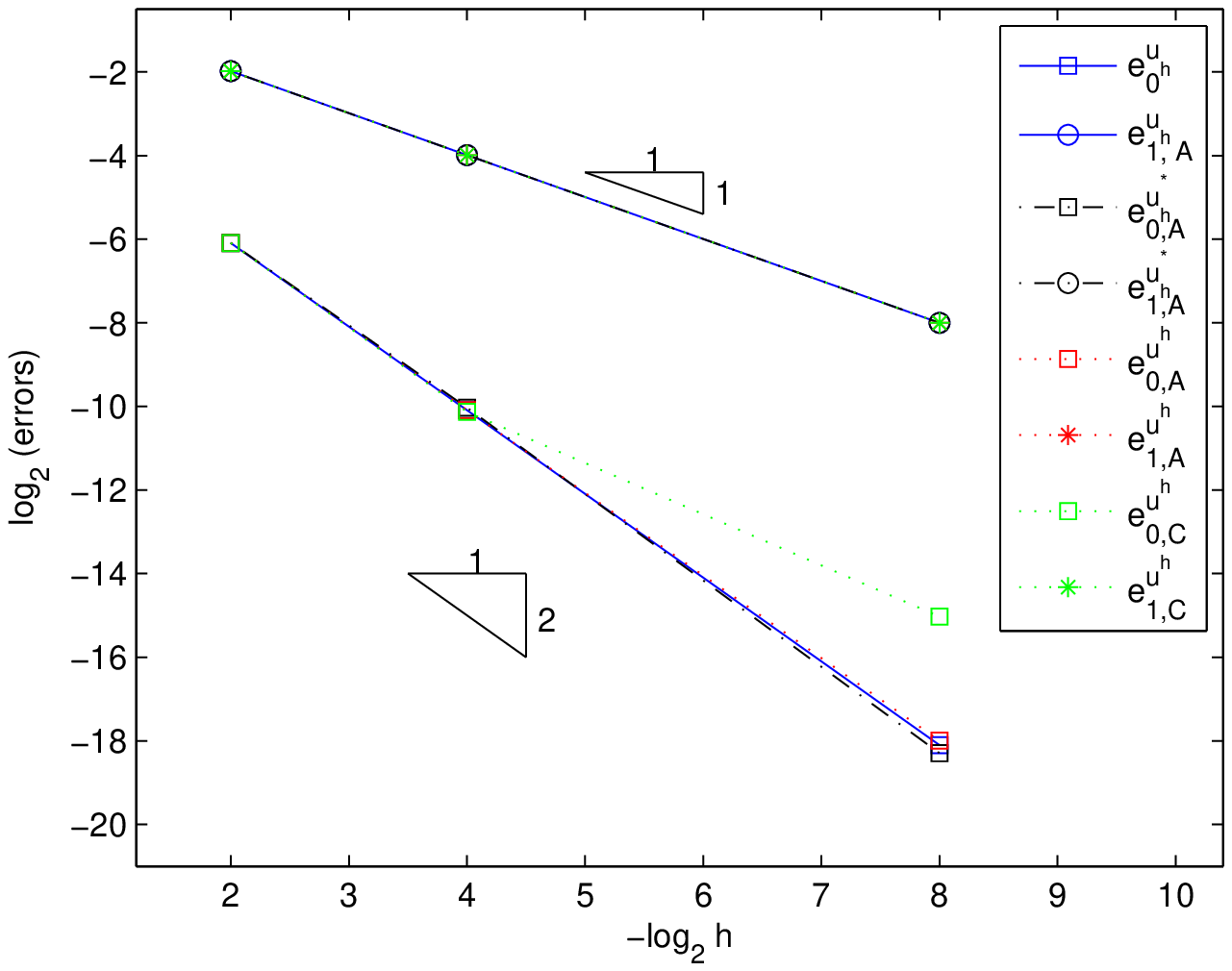}
\vskip .1cm
\includegraphics[height=5.6cm]{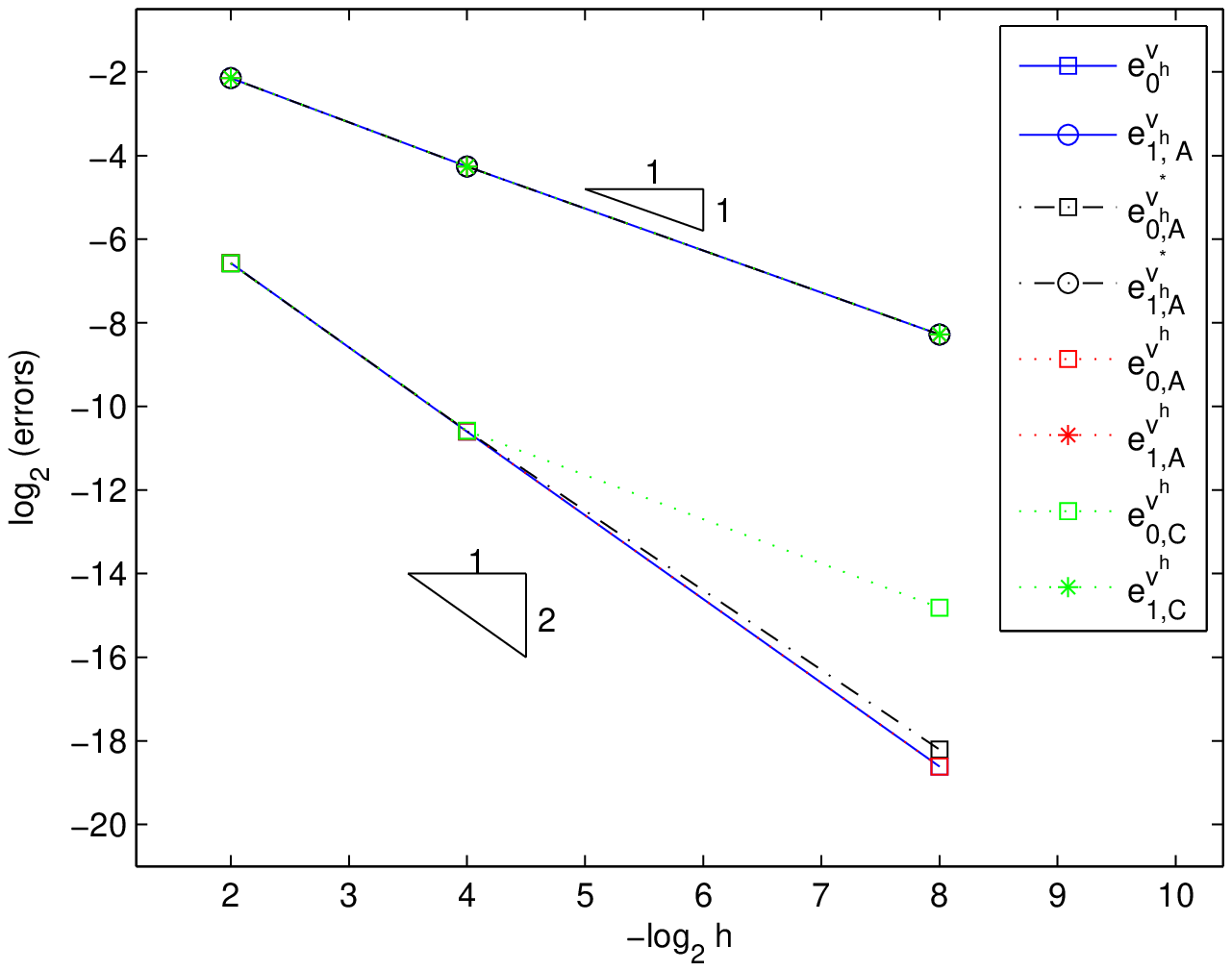}
\hskip -0.2cm
\includegraphics[height=5.6cm]{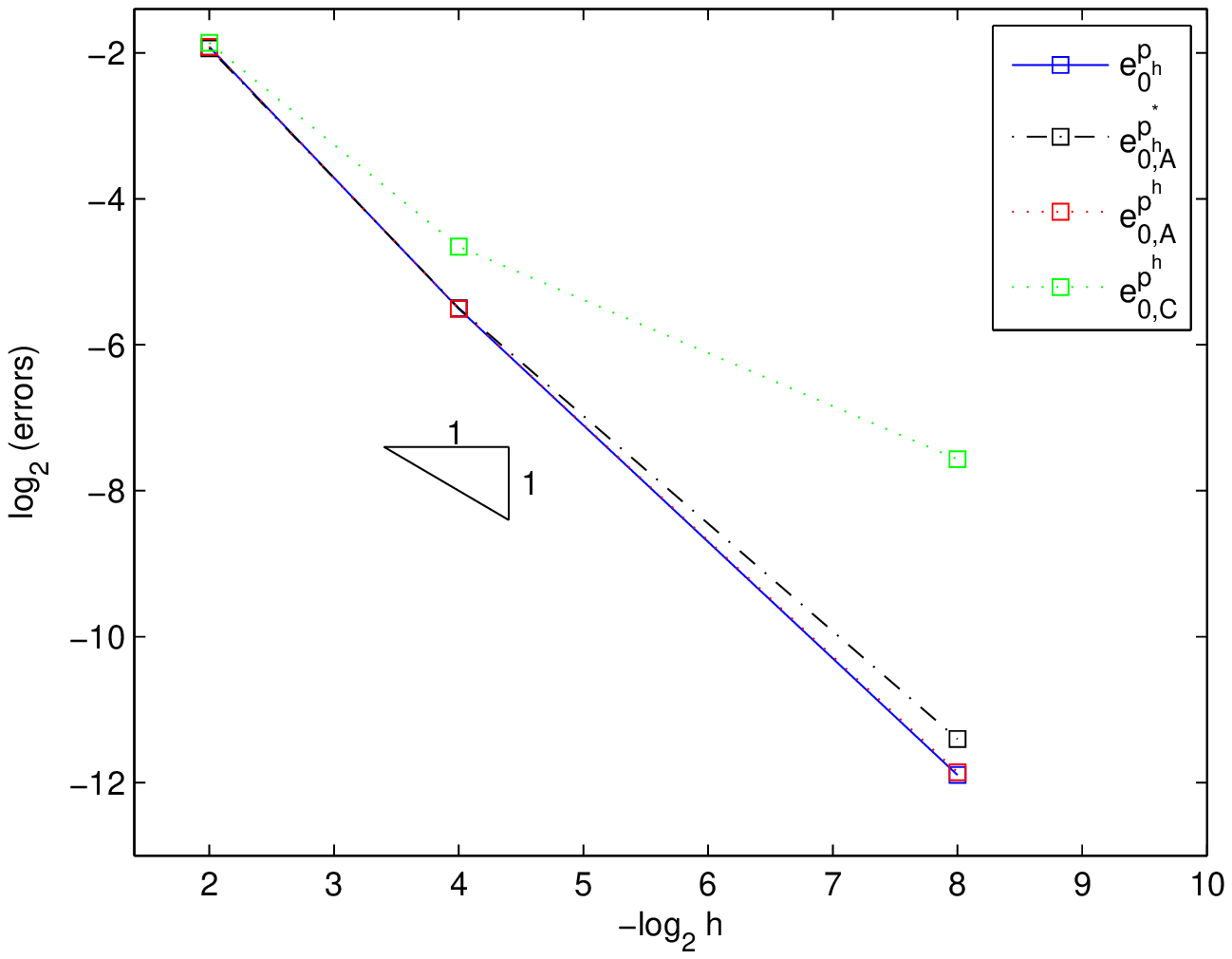}
\vskip -.1cm
\caption{Multilevel tests: plots of the errors from Table 5.3. 
Comparisons of {\it Algorithm A} and {\it Algorithm C}.
}
\label{Fig-Mini-ML_error}
\end{figure}

\begin{table}[ht]
\begin{center}
{\scriptsize
\begin{tabular}{|c||c|c||c|c||c|c||c|}
\hline
$h_l=h_{l-1}^2$   &$e^{\phi_h}_0$  &$e^{{\phi}_h}_1$ &$e^{u_h}_0$ &$e^{u_h}_1$ &$e^{v_h}_0$ &$e^{v_h}_1$ &$e^{p_h}_0$  \\
$2^{-1}$   &3.009E-3 &2.100E-2 &8.877E-3 &7.608E-2 &6.391E-3 &4.943E-2 &7.601E-2 \\
$2^{-2}$   &3.587E-4 &5.552E-3 &1.288E-3 &1.874E-2 &7.934E-4 &1.125E-2 &8.659E-3 \\
$2^{-4}$   &5.584E-6 &3.648E-4 &2.221E-5 &1.156E-3 &1.161E-5 &6.732E-4 &2.930E-4  \\
${2^{-8}}$ &1.372E-9 &1.449E-6 &5.482E-9 &4.514E-6 &2.827E-9 &2.620E-6 &1.078E-6  \\
\hline
\hline
$h_l=h_{l-1}^2$  &$e^{\phi^*_h}_{0, A}$  &$e^{\phi^*_h}_{1, A}$ &$e^{u^*_h}_{0, A}$ &$e^{u^*_h}_{1, A}$ &$e^{v^*_h}_{0, A}$ &$e^{v^*_h}_{1, A}$ &$e^{p^*_h}_{0, A}$  \\
${2^{-2}}$ &3.968E-4 &5.698E-3 &1.286E-3 &1.873E-2 &7.940E-4 &1.126E-2 &8.676E-3 \\
${2^{-4}}$ &9.394E-6 &3.743E-4 &2.219E-5 &1.156E-3 &1.165E-5 &6.732E-4 &2.932E-4 \\
${2^{-8}}$ &7.161E-8 &1.632E-6 &9.735E-9  &4.514E-6 &9.342E-9 &2.621E-6 &1.080E-6 \\
$h_l=h_{l-1}^2$  &$e^{\phi^h}_{0, A}$  &$e^{\phi^h}_{1, A}$ &$e^{u^h}_{0, A}$ &$e^{u^h}_{1, A}$ &$e^{v^h}_{0, A}$ &$e^{v^h}_{1, A}$ &$e^{p^h}_{0, A}$  \\
${2^{-2}}$ &3.585E-4 &5.552E-3 &1.288E-3 &1.874E-2 &7.934E-4 &1.125E-2 &8.659E-3 \\
${2^{-4}}$ &5.582E-6 &3.648E-4 &2.221E-5 &1.156E-3 &1.161E-5 &6.732E-4 &2.930E-4 \\
${2^{-8}}$ &6.394E-9 &1.451E-6 &5.487E-9 &4.514E-6 &2.842E-9 &2.620E-6 &1.078E-6 \\
\hline
\hline
$h_l=h_{l-1}^2$  &$e^{\phi^*_h}_{0, D}$  &$e^{\phi^*_h}_{1, D}$ &$e^{u^*_h}_{0, D}$ &$e^{u^*_h}_{1, D}$ &$e^{v^*_h}_{0, D}$ &$e^{v^*_h}_{1, D}$ &$e^{p^*_h}_{0, D}$  \\
${2^{-2}}$ &3.968E-4 &5.698E-3 &1.286E-3 &1.873E-2 &7.940E-4 &1.126E-2 &8.676E-3 \\
${2^{-4}}$ &9.394E-6 &3.743E-4 &2.219E-5 &1.156E-3 &1.165E-5 &6.732E-4 &2.932E-4 \\
${2^{-8}}$ &7.161E-8 &1.632E-6 &9.735E-9  &4.514E-6 &9.342E-9 &2.621E-6 &1.080E-6 \\
$h_l=h_{l-1}^2$  &$e^{\phi^h}_{0, D}$  &$e^{\phi^h}_{1, D}$ &$e^{u^h}_{0, D}$ &$e^{u^h}_{1, D}$ &$e^{v^h}_{0, D}$ &$e^{v^h}_{1, D}$ &$e^{p^h}_{0, D}$  \\
${2^{-2}}$ &         &         &1.286E-3 &1.873E-2 &7.941E-4 &1.125E-2 &8.681E-3 \\
${2^{-4}}$ &         &         &2.219E-5 &1.156E-3 &1.165E-5 &6.732E-4 &2.932E-4 \\
${2^{-8}}$ &         &         &9.731E-9 &4.514E-6 &9.339E-9 &2.621E-6 &1.080E-6 \\
\hline
\hline
$h_l=h_{l-1}^2$ &$e^{\phi^h}_{0, C}$  &$e^{\phi^h}_{1, C}$ &$e^{u^h}_{0, C}$ &$e^{u^h}_{1, C}$ &$e^{v^h}_{0, C}$ &$e^{v^h}_{1, C}$ &$e^{p^h}_{0, C}$  \\
${2^{-2}}$ &3.968E-4 &5.698E-3 &1.277E-3 &1.873E-2 &8.033E-4 &1.128E-2 &8.821E-3 \\
${2^{-4}}$ &8.278E-6 &3.711E-4 &2.241E-5 &1.156E-3 &1.220E-5 &6.747E-4 &2.967E-4 \\
${2^{-8}}$ &4.817E-7 &4.613E-6 &2.700E-7 &4.797E-6 &3.334E-7 &3.887E-6 &2.910E-6 \\
\hline
\end{tabular}
}
\caption{Tests for the multilevel algorithms under the Taylor-Hood/$P_2$ element discretization. Comparisons of {\it Algorithm A}, {\it Algorithm C} and {\it Algorithm D}.
}
\end{center}
\label{TH-ML-Err}
\end{table}

\begin{figure}[ht]
\centering
\includegraphics[height=5.6cm]{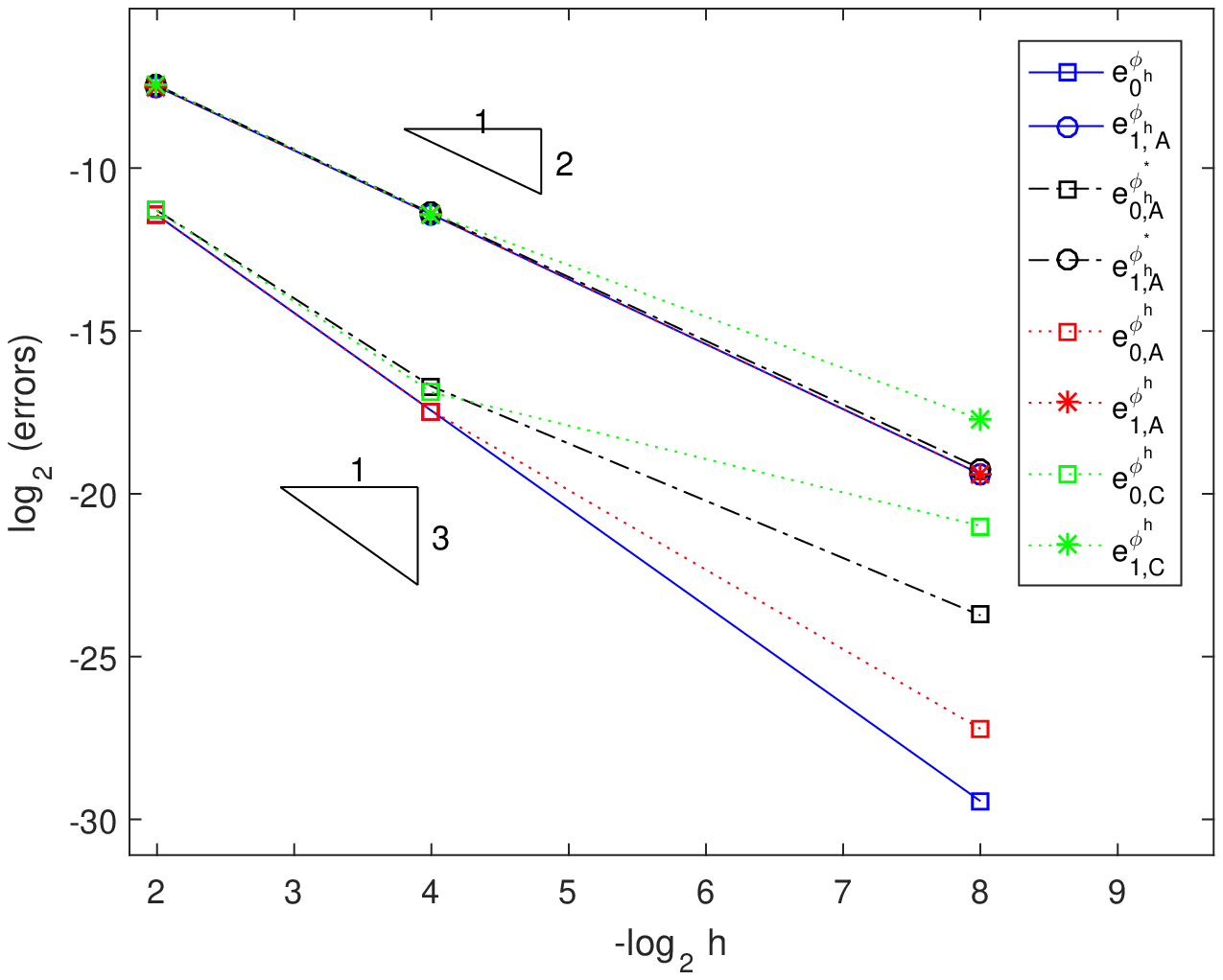}
\hskip -0.2cm
\includegraphics[height=5.6cm]{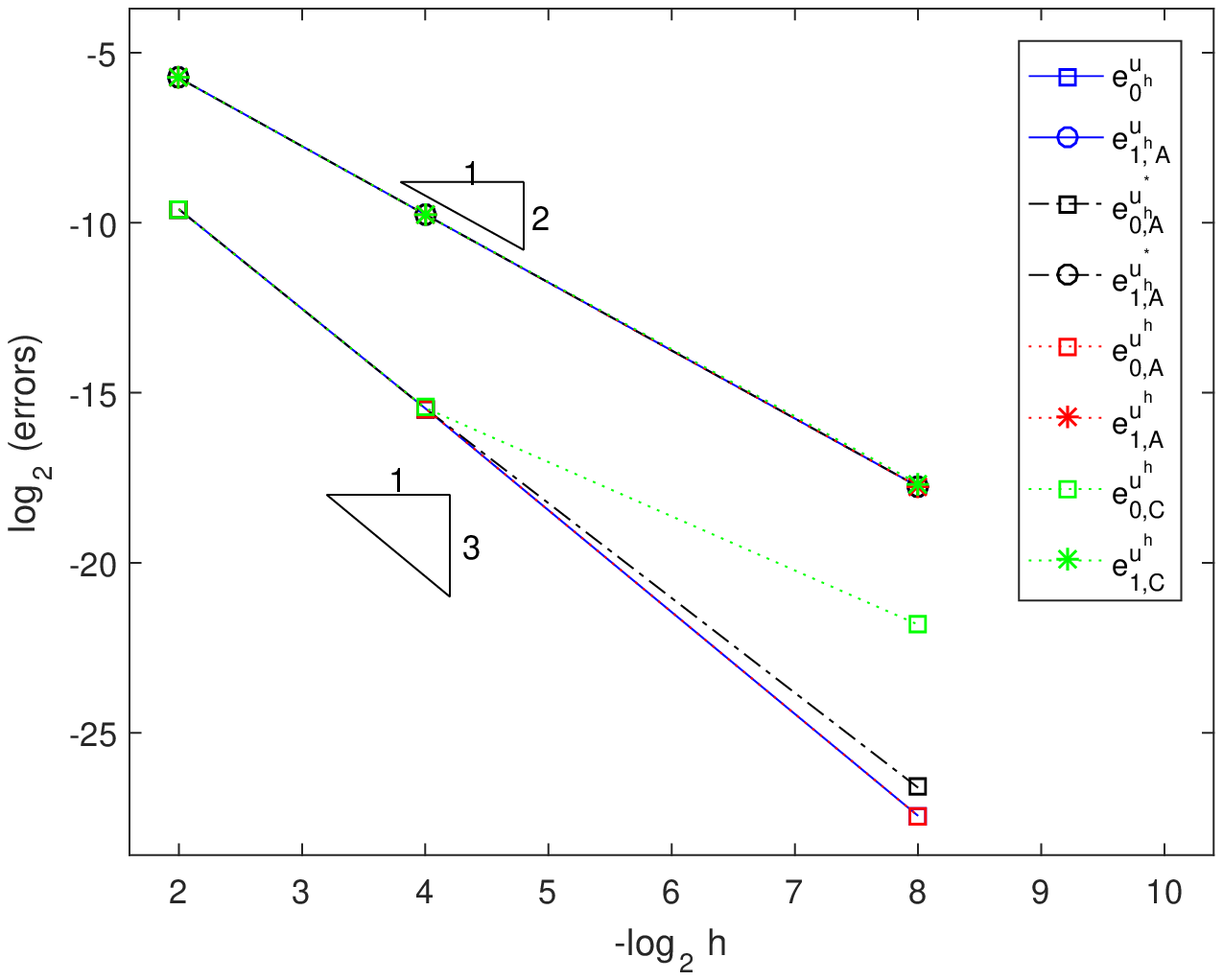}
\vskip .1cm
\includegraphics[height=5.6cm]{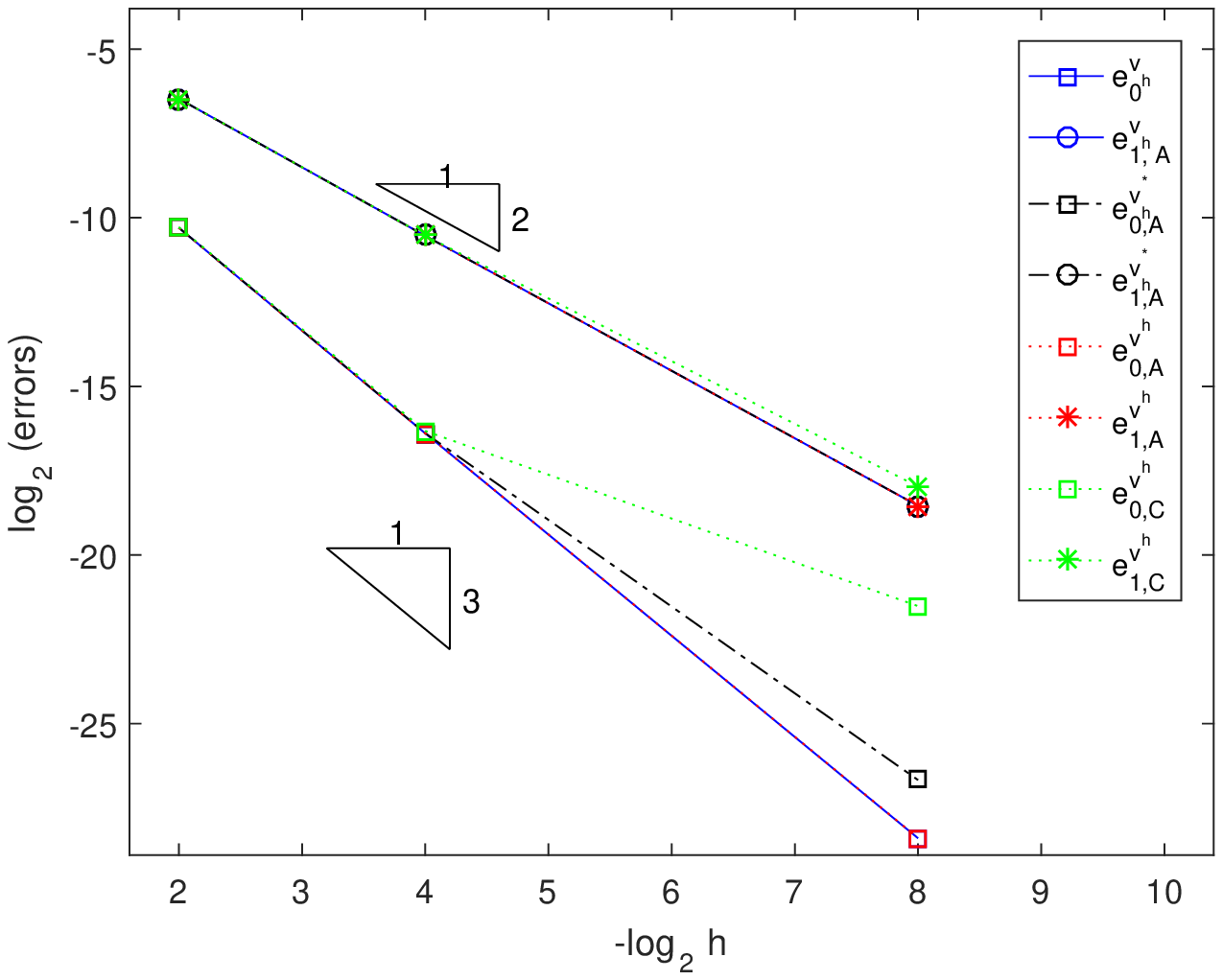}
\hskip -0.2cm
\includegraphics[height=5.6cm]{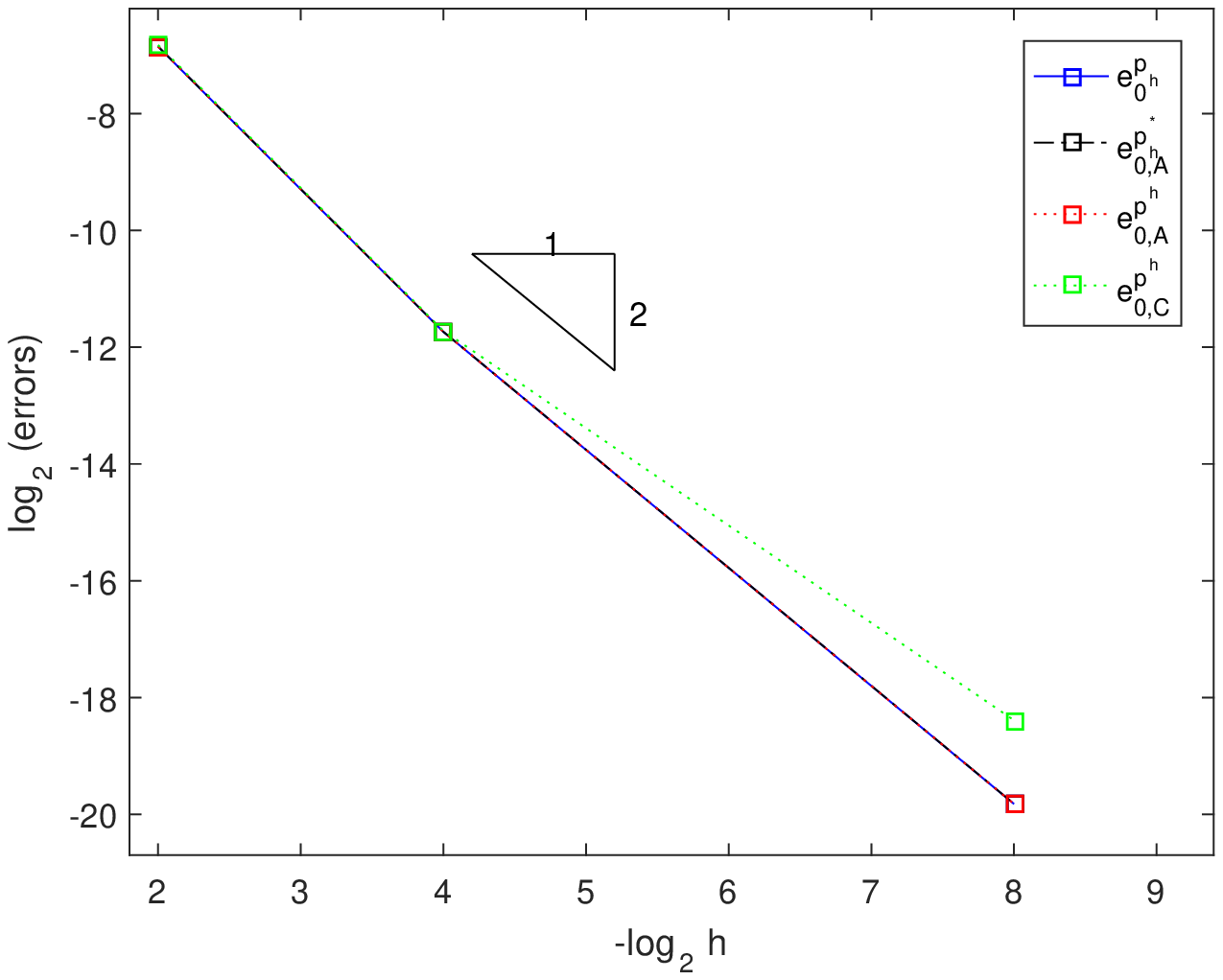}
\vskip -.1cm
\caption{Multilevel tests: plots of the errors from Table 5.4. 
Comparisons of {\it Algorithm A} and {\it Algorithm C}.
}
\label{fig-TH-ML-error}
\end{figure}

\begin{table}[htbp]
\begin{center}
{\scriptsize
\begin{tabular}{|c||c|c||c|c||c|c||c|}
\hline
$h_l=h_{l-1}^{3/2}$   &$e^{\phi_h}_0$  &$e^{{\phi}_h}_1$ &$e^{u_h}_0$ &$e^{u_h}_1$ &$e^{v_h}_0$ &$e^{v_h}_1$ &$e^{p_h}_0$  \\
$\frac{1}{2}$  &3.009E-3 &2.100E-2 &8.877E-3 &7.608E-2 &6.391E-3 &4.943E-2 &7.601E-2 \\
$\frac{1}{3}$  &8.609E-4 &9.673E-3 &2.903E-3 &3.352E-2 &1.951E-3 &2.060E-2 &2.066E-2 \\
$\frac{1}{4}$  &1.056E-4 &2.522E-3 &4.017E-4 &8.270E-3 &2.266E-4 &4.882E-3 &2.837E-3  \\
$\frac{1}{16}$ &6.774E-6 &4.146E-4 &2.692E-5 &1.315E-3 &1.410E-5 &7.663E-4 &3.361E-4 \\
$\frac{1}{58}$ &1.241E-7 &2.911E-5 &4.957E-7 &9.104E-5 &2.556E-7 &5.287E-5 &2.184E-5 \\
\hline
\hline
$h_l=h_{l-1}^{3/2}$  &$e^{\phi^*_h}_{0, D}$  &$e^{\phi^*_h}_{1, D}$ &$e^{u^*_h}_{0, D}$ &$e^{u^*_h}_{1, D}$ &$e^{v^*_h}_{0, D}$ &$e^{v^*_h}_{1, D}$ &$e^{p^*_h}_{0, D}$  \\
$\frac{1}{3}$  &9.054E-4 &9.760E-3 &2.902E-3 &3.352E-2 &1.951E-3 &2.059E-2 &2.066E-2 \\
$\frac{1}{4}$  &1.114E-4 &2.538E-3 &4.014E-4 &8.270E-3 &2.268E-4 &4.882E-3 &2.839E-3  \\
$\frac{1}{16}$ &6.848E-6 &4.150E-4 &2.692E-5 &1.315E-3 &1.411E-5 &7.663E-4 &3.361E-4 \\
$\frac{1}{58}$ &1.270E-7 &2.912E-5 &4.957E-7 &9.104E-5 &2.556E-7 &5.287E-5 &2.184E-5 \\
$h_l=h_{l-1}^{3/2}$  &$e^{\phi^h}_{0, D}$  &$e^{\phi^h}_{1, D}$ &$e^{u^h}_{0, D}$ &$e^{u^h}_{1, D}$ &$e^{v^h}_{0, D}$ &$e^{v^h}_{1, D}$ &$e^{p^h}_{0, D}$  \\
$\frac{1}{3}$  &         &         &2.901E-3 &3.352E-2 &1.952E-3 &2.059E-2 &2.069E-2 \\
$\frac{1}{4}$  &         &         &4.014E-4 &8.270E-3 &2.268E-4 &4.882E-3 &2.839E-3  \\
$\frac{1}{16}$ &         &         &2.692E-5 &1.315E-3 &1.411E-5 &7.663E-4 &3.361E-4 \\
$\frac{1}{58}$ &         &         &4.957E-7 &9.104E-5 &2.556E-7 &5.287E-5 &2.184E-5 \\
\hline
\hline
$h_l=h_{l-1}^{3/2}$ &$e^{\phi^h}_{0, C}$  &$e^{\phi^h}_{1, C}$ &$e^{u^h}_{0, C}$ &$e^{u^h}_{1, C}$ &$e^{v^h}_{0, C}$ &$e^{v^h}_{1, C}$ &$e^{p^h}_{0, C}$  \\
$\frac{1}{3}$  &9.054E-4 &9.760E-3 &2.888E-3 &3.352E-2 &1.957E-3 &2.058E-2 &2.079E-2 \\
$\frac{1}{4}$  &1.084E-4 &2.530E-3 &4.003E-4 &8.270E-3 &2.275E-4 &4.883E-3 &2.854E-3  \\
$\frac{1}{16}$ &6.814E-6 &4.149E-4 &2.691E-5 &1.315E-3 &1.411E-5 &7.663E-4 &3.364E-4 \\
$\frac{1}{58}$ &1.369E-7 &2.912E-5 &4.970E-7 &9.104E-5 &2.601E-7 &5.287E-5 &2.185E-5 \\
\hline
\end{tabular}
}
\vspace{2mm} \caption{Multilevel tests: Comparisons of {\it Algorithm A}, {\it Algorithm C} and {\it Algorithm D} using Taylor-Hood/P2 discretization. The errors between the exact solution and the solutions of the multilevel algorithms with $h_l=h_{l-1}^{3/2}$.
}
\end{center}
\label{TH-ML-Err-ADC}
\end{table}

From Table \ref{Mini-TG-Err-H3} and Table 5.2, 
we see that {\it Algorithm A} and {\it Algorithm B} actually give almost the same numerical accuracy. They only have some difference in the intermediate-step solution errors. This means that it doesn't matter whether the NS problem or the Darcy problem is solved firstly. Therefore, in the multilevel tests, we will only compare {\it Algorithm A} with {\it Algorithm C} and {\it Algorithm D}.

In Table 5.3 and Table 5.4 
we report the numerical results based on the Mini/$P_1$ element discretization and the Taylor-Hood/$P_2$ element discretization, respectively. Correspondingly, Figure 5.3 and Figure 5.4 plot the results from Table 5.3 and Table 5.4. For the both the first order and the second order discretizations, the scalings of the two successive meshsizes are all set as $h_l=h_{l-1}^2$. From Table 5.3 and Figure 5.3, we see that all the final-step solution errors based on the decoupled multilevel {\it Algorithm A} are almost the same as those based on the coupled nonlinear algorithm. This clearly shows the approximation properties of {\it Algorithm A}. For {\it Algorithm C} and {\it Algorithm D}, although they provide accurate $H^1-$ norm errors for the variables $\phi$, $u$, $v$, they can not give accurate pressure errors and the $L^2$- norm errors. To be more precisely, {\it Algorithm C} can not provide optimal $L^2$ norm errors for all variables; {\it Algorithm D} can not provide optimal $L^2$ norm errors for $\phi$.

For the second order discretization, from the results reported in Table 5.4 and Figure 5.4, we can draw the same conclusions for {\it Algorithm A} as those based on the first order discretization. For {\it Algorithm C}, both the energy norm errors and the $L^2$ norm errors are not accurate enough because the scalings for the two successive mesh sizes are set as $h_l=h_{l-1}^2$. For {\it Algorithm D}, we note that {\it Algorithm D} does provide accurate energy norm errors for fluid variables. However, it does not give the optimal errors for $\phi$ (for both the $H^1$ and $L^2$ norm errors). The reason is that one can not theoretically guarantee the optimal convergence rate of {\it Algorithm D} if $h_l=h_{l-1}^{2}$ for the second order discretization. Instead, one has to set $h_l=h_{l-1}^{3/2}$ for the second order discretization. To verify this, we report the numerical results in Table 5.5. From the results in Table 5.5, we note that both {\it Algorithm C} and {\it Algorithm D} give almost the same errors as the coupled nonlinear algorithm. The results confirm our theoretical predications, and most importantly, the results suggest that it is necessary to have the correction steps for both the Navier-Stokes subproblem and the Darcy subproblem.

\subsection{Experiments for the multilevel algorithms using different scalings between different meshlevel sizes}

\begin{table}[h]
\begin{center}
{\scriptsize
\begin{tabular}{|c||c|c||c|c||c|c||c|}
\hline
$h_l$   &$e^{\phi_h}_0$  &$e^{{\phi}_h}_1$ &$e^{u_h}_0$ &$e^{u_h}_1$ &$e^{v_h}_0$ &$e^{v_h}_1$ &$e^{p_h}_0$  \\
$2^{-1}$   &2.153E-2 &2.351E-1 &5.524E-2 &5.087E-1 &4.295E-2 &5.252E-1 &1.024E-0 \\
$2^{-3}$   &1.736E-3 &6.134E-2 &3.685E-3 &1.263E-1 &2.588E-3 &1.066E-1 &7.420E-2 \\
$2^{-6}$   &2.766E-5 &7.693E-3 &5.697E-5 &1.564E-2 &3.996E-5 &1.289E-2 &2.255E-3 \\
\hline
\hline
$h_l$  &$e^{\phi^*_h}_{0, A}$  &$e^{\phi^*_h}_{1, A}$ &$e^{u^*_h}_{0, A}$ &$e^{u^*_h}_{1, A}$ &$e^{v^*_h}_{0, A}$ &$e^{v^*_h}_{1, A}$ &$e^{p^*_h}_{0, A}$  \\
$2^{-3}$   &7.649E-3 &6.848E-2 &3.830E-3 &1.263E-1 &2.573E-3 &1.067E-1 &7.541E-2 \\
$2^{-6}$   &4.213E-4 &7.915E-3 &5.485E-5 &1.564E-3 &4.054E-5 &1.289E-2 &2.414E-3 \\   
$h_l$  &$e^{\phi^h}_{0, A}$  &$e^{\phi^h}_{1, A}$ &$e^{u^h}_{0, A}$ &$e^{u^h}_{1, A}$ &$e^{v^h}_{0, A}$ &$e^{v^h}_{1, A}$ &$e^{p^h}_{0, A}$  \\
$2^{-3}$   &1.741E-3 &6.134E-2 &3.685E-3 &1.263E-1 &2.588E-3 &1.066E-1 &7.421E-2 \\
$2^{-6}$   &2.760E-5 &7.693E-3 &5.709E-5 &1.564E-3 &3.995E-5 &1.289E-2 &2.255E-3 \\
\hline
\hline
$h_l$ &$e^{\phi^h}_{0, C}$  &$e^{\phi^h}_{1, C}$ &$e^{u^h}_{0, C}$ &$e^{u^h}_{1, C}$ &$e^{v^h}_{0, C}$ &$e^{v^h}_{1, C}$ &$e^{p^h}_{0, C}$  \\
$2^{-3}$   &7.649E-3 &6.848E-2 &3.830E-4 &1.263E-1 &2.573E-3 &1.067E-1 &7.541E-2 \\
$2^{-6}$   &4.090E-4 &7.903E-3 &1.094E-4 &1.566E-3 &1.416E-4 &1.297E-2 &2.151E-2 \\
\hline
\end{tabular}
}
\vspace{2mm} \caption{3-level test: comparisons of {\it Algorithm A} and {\it Algorithm C} using Mini/P1 discretization. The errors between the exact solution and the solutions of the multilevel algorithm with $h_1=h_{0}^{3}$ and $h_2=h_{1}^{2}$ .
}
\end{center}
\label{ML-Mini-3L}
\end{table}


\begin{table}[h]
\begin{center}
{\scriptsize
\begin{tabular}{|c||c|c||c|c||c|c||c|}
\hline
$h_l$   &$e^{\phi_h}_0$  &$e^{{\phi}_h}_1$ &$e^{u_h}_0$ &$e^{u_h}_1$ &$e^{v_h}_0$ &$e^{v_h}_1$ &$e^{p_h}_0$  \\
$1/2$   &2.153E-2 &2.351E-1 &5.524E-2 &5.087E-1 &4.295E-2 &5.252E-1 &1.024E-0 \\
$1/6$   &1.056E-4 &2.522E-3 &4.017E-4 &8.270E-3 &2.266E-4 &4.882E-3 &2.837E-2 \\
$1/36$  &4.918E-7 &7.277E-5 &1.965E-6 &2.282E-4 &1.015E-6 &1.326E-4 &5.515E-5 \\
\hline
\hline
$h_l$  &$e^{\phi^*_h}_{0, A}$  &$e^{\phi^*_h}_{1, A}$ &$e^{u^*_h}_{0, A}$ &$e^{u^*_h}_{1, A}$ &$e^{v^*_h}_{0, A}$ &$e^{v^*_h}_{1, A}$ &$e^{p^*_h}_{0, A}$  \\
$1/6$   &1.808E-4 &2.830E-3 &4.006E-4 &8.270E-3 &2.275E-4 &4.883E-3 &2.849E-3 \\
$1/36$  &1.300E-6 &7.463E-5 &1.966E-6 &2.282E-4 &1.017E-6 &1.326E-4 &5.521E-5 \\
$h_l$  &$e^{\phi^h}_{0, A}$  &$e^{\phi^h}_{1, A}$ &$e^{u^h}_{0, A}$ &$e^{u^h}_{1, A}$ &$e^{v^h}_{0, A}$ &$e^{v^h}_{1, A}$ &$e^{p^h}_{0, A}$  \\
$1/6$   &1.055E-4 &2.522E-3 &4.018E-4 &8.270E-3 &2.266E-4 &4.882E-3 &2.837E-3 \\
$1/36$  &4.921E-7 &7.277E-5 &1.965E-6 &2.282E-4 &1.015E-6 &1.326E-4 &5.515E-5 \\
\hline
\hline
$h_l$ &$e^{\phi^h}_{0, C}$  &$e^{\phi^h}_{1, C}$ &$e^{u^h}_{0, C}$ &$e^{u^h}_{1, C}$ &$e^{v^h}_{0, C}$ &$e^{v^h}_{1, C}$ &$e^{p^h}_{0, C}$  \\
$1/6$   &1.808E-4 &2.830E-3 &4.006E-4 &8.270E-3 &2.275E-4 &4.883E-3 &2.849E-3 \\
$1/36$  &1.121E-5 &1.225E-4 &8.057E-6 &2.310E-4 &8.614E-6 &1.470E-4 &8.067E-5 \\
\hline
\end{tabular}
}
\vspace{2mm} \caption{3-level tests: comparisons of {\it Algorithm A} and {\it Algorithm C} using Taylor-Hood/$P_2$ discretization. The errors between the exact solution and the solutions of the multilevel algorithm with $h_1 \approx h_{0}^{5/2}$ and $h_2=h_{1}^{2}$ .
}
\end{center}
\label{ML-TH-3L}
\end{table}


From \cite{huang2016newton}, we see that for the first two levels of the multilevel {\it Algorithm A}, one can take $h_1=h_0^3$ to guarantee the optimal convergence of the energy norm errors (for simplicity, we use the first order discretization for the discussion). However, the analysis in Section 4 shows that one should take $h_l=h_{l-1}^2$ to guarantee the solution errors are optimal in the energy norm on all mesh levels. This suggests us to test the multilevel algorithms using different scalings on different meshlevels. In this subsection, we test the multilevel algorithms under the three level cases using different scalings between $h_l$ and $h_{l-1}$ for two adjacent mesh levels.

For the first order discretization, we take $h_1=h_0^3$ while $h_2=h_1^2$. The corresponding numerical results are reported in Table 5.6. For the second order discretization, we set $h_1=h_0^{5/2}$ while $h_2=h_1^2$. The corresponding numerical results are reported in Table 5.7. 
From both Table 5.6 and Table 5.7, by comparing the FE errors with the multilevel algorithm errors, we see that {\it Algorithm A} still gives optimal energy norm errors and optimal $L^2$- norm errors for velocity. For {\it Algorithm C}, as the scaling between $h_l$ and $h_{l-1}$ is of very higher order, neither theoretical analysis nor numerical experiments guarantee it can give optimal energy norm or $L^2$ norm solution errors.

\section{Conclusion}
In conclusion, we have proposed some decoupled and linearized multilevel algorithms for the coupled NS/Darcy model. These algorithms are numerically efficient and also enables easy and efficient implementation and software reuse. Numerical analysis are presented to show that the decoupled and linearized {\it Algorithm A} retains the same order of approximation accuracy as the coupled and nonlinear algorithm if the scalings between two successive mesh level sizes are properly chosen. Extensive numerical experiments are provided to verify the theoretical predictions and to compare the different algorithms.


\penalty-8000

%
%

\end{document}